\documentclass[11pt,reqno]{amsart}
\usepackage{graphicx}
\usepackage{amssymb,latexsym,epsfig,mathrsfs,graphpap,graphics,amsmath,amssymb}
\usepackage{amstext}
\usepackage{verbatim}
\usepackage{color}

\newtheorem{conj.}[thm]{Conjecture}

\theoremstyle{definition}

\theoremstyle{remark}

\numberwithin{equation}{section}

\begin{document}

\begin{flushleft}
  {\bf\Large {Nonuniform  Multiresolution Analysis Associated with\\[1.8mm]  Linear Canonical Transform}}
\end{flushleft}

\parindent=0mm \vspace{.3in}

  {\bf{ Firdous A. Shah$^{\star}$ and  Waseem Z. Lone$^{\star}$  }}

\parindent=0mm \vspace{.1in}
{\small \it $^{\star\star}$Department of  Mathematics, University of Kashmir, South Campus,
Anantnag 192101, Jammu and Kashmir, India.
E-Mail: $\text{fashah@uok.edu.in}$;\,$\text{lwaseem.scholar@kashmiruniversity.net}$.}

\parindent=0mm \vspace{.1in}
{\small {\bf Abstract.} The linear canonical transform (LCT) has attained respectable status within a short span and is being broadly employed across several disciplines of science and engineering including signal processing, optical and radar systems, electrical and communication systems, quantum physics etc, mainly due to the extra degrees of freedom and simple geometrical manifestation. In this article, we introduce a novel multiresolution analysis (LCT-NUMRA) on the spectrum $\Omega = \big\{0,{r}/{N}\big\}+2\mathbb{Z},$ where  $N \geqq 1$ is an integer and $r$ is an odd integer with $ 1 \leqq r \leqq 2N-1,$ such that $r$ and $N$ are relatively prime, by  intertwining the ideas of nonuniform MRA and linear canonical transforms. We first develop nonuniform multiresolution analysis associated with the LCT and then drive an algorithm to construct an LCT-NUMRA starting from a linear canonical low-pass filter $\Lambda_{0}^{M}(\omega)$ with suitable conditions. Nevertheless, to extend the scope of the present study, we construct the associated wavelet packets for such an MRA and investigate their properties by means of the linear canonical transforms.  

\parindent=0mm \vspace{.1in}

\parindent=0mm \vspace{.1in}
 }

\section{Introduction}

\parindent=0mm \vspace{.0in}

In the early 1970s, a promising linear integral transform with three free parameters, namely, linear canonical transform was independently introduced by Collins \cite{C} in paraxial optics, and Moshinsky, and Quesne \cite{MQ} in quantum mechanics, to study the conservation of information and uncertainty under linear maps of phase space. The LCT  provides a unified treatment of the generalized Fourier transforms in the sense that it is an embodiment of several well-known integral transforms including the Fourier transform, fractional Fourier transform, Fresnel transform,  scaling operations and so on \cite{A,JA,TDW}. Over a couple of decades, the application areas for LCT have been growing at an exponential rate and is as such befitting for investigating deep problems in time-frequency analysis, filter design, phase retrieval problems, pattern recognition, radar analysis, holographic three-dimensional television, quantum physics, and many more. Apart from applications, the theoretical skeleton of LCT has likewise been extensively studied and investigated including the convolution theorems, sampling theorems, Poisson summation formulae, uncertainty principles,  shift-invariant theory and so on. For more about LCT and their applications, we allude to  \cite{BS,HKOS,XL}.

\parindent=8mm \vspace{.1in}
The celebrated Mallat's theory of multiresolution (MRA)  is one of the most remarkable, significant and elegant  tool for constructing orthonormal wavelet basis for $L^2(\mathbb R)$ \cite{M}. This remarkable idea brings with it a new way of thinking, which was entirely missing in previously existing algorithms and has been the primary source of many new evolutions in wavelet analysis.  Some popular wavelets such as Shannon wavelet, Meyer wavelet,  Franklin wavelet, spline wavelets, biorthogonal wavelets, harmonic wavelets and Daubechies wavelets have been constructed by employing the notion of an MRA \cite{D}. It is a well-known fact that a new wavelet is born each day to take the place of the one that does die on any given day and as a consequence, some off-shoots of the MRA have been formulated for constructing different types of wavelet bases in recent years, for instance, frame MRA, generalized MRA, multiscale MRA, periodic MRA, non-stationary MRA, vector-valued MRA, two-direction MRA, $p$-adic MRA, nonuniform MRA and many others \cite{DS1,DS2}. Keeping in view the fact that all  these MRA's are developed on regular lattices, Gabardo and Nashed \cite{G1} presented an unprecedented MRA for the dilation $2N$ and the translation set $\Omega$ given by $\Omega=\left\{0, \frac{r}{N}\right\}+2 \mathbb{Z},$ where $N\geqq 1$ is an integer and $r$ is an odd integer with $1\leqq r \leqq 2N-1,$ such that $r$ and $N$ are relatively prime, acting on the scaling function $\phi$, is no longer a group, but a union of two lattices, which is associated with a famous Fuglede conjecture on spectral pairs.  A nonuniform MRA (NUMRA) is a non-decreasing collection of closed subspaces $\left\{ V_j: j\in\mathbb Z\right\}$ of $L^2(\mathbb R)$ satisfying the following axioms:

\parindent=0mm \vspace{.0in}
(i)\; $V_j\subset V_{j+1},\quad j\in \mathbb{Z}$;

\parindent=0mm \vspace{.0in}
(ii)\; $\bigcup_{j\in\mathbb Z} V_j$ is dense in $L^2 (\mathbb{R})$  and $\bigcap_{j\in\mathbb{Z}} V_j=\left\{0\right\}$;

\parindent=0mm \vspace{.0in}
(iii)\; $f(\cdot)\in V_j$ if and only if $f(2N\cdot)\in V_{j+1}$;

\parindent=0mm \vspace{.0in}
(iv) There exists a function $\phi\in V_0$ such that $\left\{\phi(\cdot-\lambda):\lambda\in \Omega \right\}$ is an orthonormal basis for $ V_0$.

\parindent=0mm \vspace{.1in}
It is pertinent to note that when $N=1$ and $\Omega =\mathbb {Z}$, the nested family $\left\{V_j:j\in\mathbb Z\right\}$   boils down to the classical definition of an MRA with dyadic dilation. These studies were proceeded by Gabardo and his colleagues in \cite{G2,G3,G4}, wherein they establish an extension of Cohen's theorem which provides the necessary and sufficient condition for the orthonormality of the system $\left\{\phi(\cdot-\lambda):\lambda\in \Omega \right\}$ and presented some equivalent conditions of the associated wavelets via dimension functions. The theory of nonuniform wavelets was further studied and extensively investigated by many researchers in different directions, for example, nonuniform wavelet packets \cite{Beh}, vector-valued nonuniform wavelet packets \cite{MM}, generalized nonuniform MRA \cite{MN}, nonuniform wavelet frames  \cite{Shah1,VP}, nonuniform wavelets and wavelet packets on local fields of positive characteristic \cite{Shah3,Shah4,Shah5}.

\pagestyle{myheadings}

\parindent=8mm \vspace{.1in}
Keeping in view the fact that  LCT is endowed with three degrees of freedom, which in turn enhance the signal processing capability to non-orthogonal phase-space directions, therefore it is intriguing to present a novel multiresolution analysis in $L^2(\mathbb R)$ with the aim of constructing a new class of nonuniform wavelets which are comparatively flexible and enjoys additional degrees of freedom and consequently,  provides a richer representation of non-stationary signals. In this setup,
the associated core subspace $V_{0}$ of $L^2(\mathbb R)$ has an orthonormal basis of the form  $\big\{ \phi( t - \lambda ) \, e^{ -\frac{i \pi A}{B} ( t^2 - \lambda^2 ) }: \lambda \in \Omega, M=(A,B,C,D),B\neq 0\big \}$, where the translation set $\Omega$ is not necessarily a group, but it is the union of $\mathbb{Z}$ and a translate of $\mathbb{Z}$. To facilitate the motive, we establish a necessary and sufficient condition for the existence of the  linear canonical nonuniform wavelets  and develop an algorithm to construct an LCT-NUMRA starting from a linear canonical low-pass filter $\Lambda_{0}^{M}(\omega)$ with suitable conditions. To broaden the scope of the present study, we introduce the notion of nonuniform wavelet packets associated with the linear canonical transforms by adopting the procedure of splitting the wavelet spaces pioneered by Coifman et al.\cite{CM}. From the outset, we investigate the orthogonal properties of the associated wavelet packets using the machinery of LCT and, then we construct several new orthonormal bases for $L^2(\mathbb R)$ by means of the orthogonal decomposition of wavelet packet subspaces.  

\parindent=8mm \vspace{.1in}
The rest of the article is structured as follows: In Section 2, we first recapitulate the linear canonical transform and then introduce the novel MRA associated with LCT on the spectrum $\Omega$. Section 3 is entirely devoted to the construction of an LCT-NUMRA starting from a linear canonical low-pass filter $\Lambda_{0}^{M}(\omega)$. In Section 4, we construct linear canonical wavelet packets and study their orthogonal properties via LCT. Finally, in Section 5, we present the orthogonal decomposition of wavelet packet subspaces to construct many orthonormal bases for $L^2(\mathbb R)$.

\section{ Nonuniform  Multiresolution Analysis Associated with LCT }

\parindent=0mm \vspace{.0in}
We shall start this section with a brief overview of the linear canonical transform  which serves as a cornerstone for the subsequent developments and then introduce the notion of the LCT-NUMRA in the context of time-frequency analysis. Subsequently, we provide a complete characterization of the associated orthonormal nonuniform wavelets in terms of the linear canonical low-pass filter.

\parindent=0mm \vspace{.0in}
{\it 2.1. The Linear Canonical Transform}

\parindent=0mm \vspace{.0in}
For the sake of brevity, we may write the matrix $M_{2\times 2}$ as $M=\left(A,B,C,D\right)$ and the corresponding transpose matrix shall be denoted as $M^{T}=\left(A,B,C,D\right)^{T}$.

\parindent=0mm \vspace{.1in}
{\bf Definition 2.1.} The linear canonical transform of any $f\in L^2(\mathbb R)$ with respect to the uni-modular matrix $M=(A,B,C,D)$  is defined by
\begin{align*}
\mathscr L_A\big[f\big](\omega)=\left\{\begin{array}{cc}\displaystyle\int_{\mathbb R}f(x)\,{\mathcal K}_M(t,\omega)\,dt,\,& B\neq 0\\
 \sqrt{D}\, \exp\left\{\dfrac{CD\omega^2}{2}\right\} f(D\omega ),& B= 0,\\ \end{array}\right.\tag{2.1}
\end{align*}
where $K_M(t,\omega)$ is called the kernel of the LCT  and is given by
\begin{align*}
{\mathcal K}_M(t,\omega)=\dfrac{1}{\sqrt{2i \pi B}}\,\exp\left\{\dfrac {i(At^2-2t\omega+D\omega^2)}{2B}\right\},\quad B\neq 0.\tag{2.2}
\end{align*}

\parindent=0mm \vspace{.0in}
It is worth noticing that for the case $B=0$, the LCT (2.1) corresponds to a chirp multiplication operation and the same is of no particular interest to us. As such, in the rest of the article, we shall focus our attention on the case $B\neq 0$. We also note that the phase-space transform (2.1) is lossless if and only if the matrix $M$ is unimodular; that is, $AD-BC=1$ \cite{XL}. The linear canonical transform (2.1) includes several known transforms as special cases. For example, for $M=(1,B,0,1)$, we obtain the Fresnel transform, for $M=(\cos\theta, \sin\theta, -\sin\theta, \cos\theta,0,0)$, the LCT  boils down to the fractional Fourier transform whereas for $M=(0,1,-1,0,0,0)$, we obtain the classical Fourier transform. Moreover,  Bi-lateral Laplace,  Gauss-Weierstrass, and Bargmann transform are also special cases of LCT \cite{BS}. Nevertheless, the inversion formula corresponding to (2.1) is given by
\begin{align*}\label{2.3}
f(t)=\int_{\mathbb R}\mathscr L_M\big[f\big](\omega)\,\overline{\mathcal{K}_M(t,\omega)}\,d\omega.\tag{2.3}
\end{align*}
Also, the Parseval's formula associated with (2.1) reads:
\begin{align*}
\Big\langle \mathscr L_M\big[f\big], \mathscr L_M\big[g\big] \Big\rangle =\Big\langle f, g\Big\rangle,\quad \forall~ f,g\in L^2(\mathbb R).\tag{2.4}
\end{align*}

\parindent=0mm \vspace{.1in}
{\it 2.2. A Novel MRA Associated with Linear Canonical Transform}

\parindent=0mm \vspace{.1in}
For an integer $N\ge 1$ and an odd integer $r$ with $1 \le r \le 2N-1$ such that $r$ and $N$ are relatively prime, we define
\begin{align*}
\Omega=\left\{0, \dfrac{r}{N}\right\}+2\mathbb Z=\left\{ \dfrac{rk}{N}+2n: n\in \mathbb Z, k=0,1\right\}.\tag{2.5}
\end{align*}

\parindent=0mm \vspace{.0in}
It is easy to verify that $\Omega$ is neither a group nor a uniform discrete set, but is the union of $\mathbb Z$ and a translate of $\mathbb Z$. In deed, $\Omega$ is the spectrum for the spectral set $\Delta=\big[0, \frac{1}{2}\big)\cup \big[\frac{N}{2}, \frac{N+1}{2}\big)$ and the pair $(\Omega,\Delta)$ is called a  spectral pair \cite{G1}.

\parindent=8mm \vspace{.1in}
Next, we shall formally  introduce a novel  multiresolution analysis  associated with linear canonical transform on the spectrum $\Omega$ as follows:

\parindent=0mm \vspace{.1in}
{\bf{Definition 2.2.}} Given a real uni-modular matrix $M=(A,B,C,D)$ and integers $N \ge 1$ and $r$ odd  with $1\leq r \leq 2N-1$ such that $r$ and $N$ are relatively prime, an associated  linear canonical nonuniform multiresolution analysis (abbreviated LCT-NUMRA) is a collection  $\big\{ V_{j}^{M}:j\in\mathbb Z\big\}$ of closed subspaces of $L^2(\mathbb R)$ satisfying the following properties:

\parindent=0mm \vspace{.1in}
(a)\quad $V_j^M \subset V_{j+1}^M\; \text{for all}\; j \in \mathbb Z;$

\parindent=0mm \vspace{.1in}
(b)\quad $\bigcup_{j\in \mathbb Z}V_j^M\;\text{is dense in}\;L^2({\mathbb R});$

\parindent=0mm \vspace{.1in}
(c)\quad $\bigcap_{j\in \mathbb Z}V_j^M=\{0\};$

\parindent=0mm \vspace{.0in}
(d)\quad $f(t) \in V_j^M\; \text{if and only if}\;f(2Nt)\, e^{-i \pi A\big(1-(2N)^2\big)t^2/B} \in V_{j+1}^M\; \text{for all}\; j \in \mathbb Z;$

\parindent=0mm\vspace{.0in}

(e)~ There exists a function $\phi$ in $V_0^M$ such that $\big\{\phi^M_{ 0,\lambda }(t) = \phi( t - \lambda ) \, e^{ -\frac{i \pi A}{B} ( t^2 - \lambda^2 ) }: \lambda \in \Omega\big \}$, is a complete orthonormal basis for $V_0^M$.

\parindent=0mm\vspace{.1in}
It is worth noting that Definition 2.2 embodies certain existing MRA's and also give birth to some novel multiresolution analyses which are yet to be reported in
the open literature:

\begin{itemize}

\item For the parameter matrix $M=(\cos\theta, \sin\theta, -\sin\theta, \cos\theta )$, $\theta\neq n\pi$, we can obtain a novel fractional nonuniform MRA.

\item For the matrix $M=(1,B,0,1)$, we can obtain a novel Fresnel nonuniform MRA.

\item As a special case, when $M=(0,1,-1,0)$, the Definition 2.2 boils down to standard nonuniform MRA introduced by Gabardo and Nashed \cite{G1}.

\item For $ N = 1 $ and $M = (0,1,-1,0) $, the Definition 2.2 reduces to the classical multiresolution analysis with dilation factor equal to $2$ \cite{M}.

\item For $ N = 1 $ and $M =(\cos\theta, \sin\theta, -\sin\theta, \cos\theta )$, the Definition 2.2 reduces to the fractional MRA with dyadic dilations \cite{SLZ}.

\item For $ N = 1 $ and $M =(A,B,C,D)$, the Definition 2.2 reduces to the  MRA  associated with linear canonical transform recently introduced by \cite{WWW}.

\end{itemize}

\parindent=0mm \vspace{.1in}
Since $\phi\in V_{0}^{M} \subset V_{1}^{M}$ and the collection $\big\{\phi_{1,\lambda}^{M}:\lambda\in \Omega\big\}$ is an orthonormal basis in $V_{1}^{M}$, hence, the function $\phi\in V_{1}^{M}$  has the Fourier expansion as
\begin{align*}
\phi (t) = \sum_{ \lambda \in \Omega} a_\lambda \, \phi^M_{ 1 , \lambda }(t)  =\sqrt{2N} \sum_{ \lambda \in \Omega} a_\lambda  \phi \big( 2Nt - \lambda \big) \, e^{- \frac{i \pi A}{B}( t^2 - \lambda^2 ) }, \tag{2.6}
\end{align*}
where
\begin{align*}
 a_\lambda = \int_{\mathbb R} \phi(t) \, e^{ -i \pi A t^2/{B} } \, \overline{ \phi^M_{ 1 , \lambda }(t) } \, dt \quad \text{ and } \quad \sum_{ \lambda \in \Omega } | a_\lambda |^2 < \infty . \tag{2.7}
\end{align*}

Implementing the linear canonical transform on both sides of (2.6), we obtain
\begin{align*}
{\mathscr L}_{M}\big[\phi(t)\big](\omega)= \widehat{ \phi } \left( \frac{\omega}{B} \right) = \Lambda^M_0 \left( \frac{\omega}{2NB} \right) \widehat{ \phi } \left( \frac{\omega}{2NB} \right),\tag{2.8}																				\end{align*}
where
\begin{align*}
 \Lambda^M_0 \left( \frac{\omega}{B} \right) = \dfrac{1}{\sqrt {2N}}\sum_{ \lambda \in \Omega } a^M_\lambda \, e^{ -2 \pi i \lambda {\omega}/{B} }.\tag{2.9}    																						\end{align*}
 In view of the specific form of $\Omega= \left\{0, r/N\right\}+2\mathbb Z$, the low-pass filter  $\Lambda^M_0 $  associated with the scaling function $\phi$ can be written as
\begin{align*}
 \Lambda^M_0 \left( \frac{\omega}{B} \right) = \Lambda^{M , 1}_0 \left( \frac{\omega}{B} \right) + e^{ - 2 \pi i {\omega r}/{NB}} \, \Lambda^{ M , 2 }_0 \left( \frac{\omega}{B}\right),  \tag{2.10}
\end{align*}
where $\Lambda^{M,1 }_0 $ and $\Lambda^{M,2 }_0 $ are locally $L^2$, $B/2$-periodic functions associated with uni-modular matrix $M$. Replacing ${\omega}/{B} $ by ${\omega}/{2NB} $ in relation (2.9), we obtain
\begin{align*}
 \widehat{ \phi } \left(\frac{\omega}{2NB} \right) =\Lambda^M_0 \left(\frac{\omega}{(2N)^2 B}\right) \widehat{ \phi } \left( \frac{\omega}{(2N)^2 B}\right),
\end{align*}
and then
\begin{align*}
 \widehat{ \phi }\left( \frac{\omega}{B} \right) = \Lambda^M_0 \left( \frac{\omega}{2NB} \right) \Lambda^M_0 \left( \frac{\omega}{(2N)^2 B}  \right) \widehat{ \phi } \left( \frac{\omega}{(2N)^2 B}  \right).
\end{align*}
Continuing like this, for any $n\in\mathbb N$, we can obtain
\begin{align*}
 \widehat{ \phi }\left( \frac{\omega}{B} \right) = \widehat{ \phi } \left( \frac{\omega}{(2N)^n B}  \right) \prod_{j=1}^{n} \Lambda^M_0 \left( \frac{\omega}{(2N)^j B}  \right).
\end{align*}

\parindent=0mm\vspace{.0in}
This suggests that we should be able to let $n\rightarrow\infty$ ($\hat \phi$ to be continuous at the origin) to obtain
\begin{align*}
 \widehat{ \phi } \left( \frac{\omega}{B} \right) = \widehat{ \phi }( 0 ) \prod_{j= 1 }^{ \infty } \Lambda^M_0 \left( \frac{\omega}{(2N)^j B}  \right),\tag{2.11}       																			\end{align*}
since $ (2N)^{-j} \rightarrow 0 $ as $j\rightarrow \infty $.

\parindent=8mm\vspace{.1in}
If $ \widehat{ \phi }( 0 ) = 0 $, then $ \widehat{ \phi } \big({\omega}/{B} \big) = 0 $ for all $\omega$, and hence $\phi$ is the zero function. Therefore, a non trivial solution must satisfy $ \widehat{ \phi } ( 0 ) \neq 0 $. Assume that $ \widehat{ \phi } (0) = 1 $. Then equation (2.11) becomes
\begin{align*}
 \widehat{ \phi } \left( \frac{\omega}{B}\right) = \prod_{j=1}^{\infty} \Lambda^M_0 \left( \frac{\omega}{(2N)^j B}  \right)\tag{2.12}
\end{align*}

\parindent=0mm\vspace{.0in}
Since $ \widehat { \phi } ( 0 ) = 1 $, it follows immediately from equation (2.8) that $\Lambda^M_0(0)=1$, which is essential for convergence of the infinite product $ \prod_{j=1}^{\infty} \Lambda^M_0 \left( \frac{\omega}{(2N)^j B}  \right)  $.

\parindent=8mm \vspace{.1in}
For each $j\in\mathbb Z$ and real matrix $M=(A,B,C,D)$, the LCT wavelet subspace $W_{j}^{M}$ is defined as the orthogonal complement of $V_{j}^{M}$ in $V_{j+1}^{M}$, so that $W_{j}^{M}\perp V_{j}^{M}$. It is clear from the conditions (a), (b) and (c) of the Definition 2.2 that
\begin{align*}
 L^2(\mathbb R) = \bigoplus_{j\in \mathbb Z} W^M_{j}. \tag{2.13}
\end{align*}

\parindent=0mm \vspace{.0in}
{\bf Definition 2.3.}  A collection of square integrable functions $ \left \{ \psi^M_k :k=1,2,\dots, 2N-1  \right \} $ in $V_{1}^{M}$ will be called a set of wavelets associated with given LCT-NUMRA if the family of functions $ \big \{ \psi_k ( t - \lambda ) \, e^{ -i \pi \frac{A}{B}( t^2 - \lambda^2 ) }: 1 \leq k \leq 2N-1, \, \lambda \in \Omega\big \} $ forms an orthonormal basis for $ W^M_0 $.

\parindent=8mm \vspace{.1in}
Assume that there exists $(2N-1)$ functions $\big\{ \psi_{1}^{M}, \psi_{2}^{M},\dots, \psi_{2N-1}^{M}\big\}$ in $L^2(\mathbb R)$ such that their translates by the elements of $\Omega$ and dilations by the integer powers of $2N$ form a Riesz basis of $W_{j}^{M}$, i.e.,
\begin{align*}
W_{j}^{M}=\overline{\text{Span}}\left\{\psi_{k,j,\lambda}^{M}: 1 \leq k \leq 2N-1, \, \lambda \in \Omega \right\},~j\in\mathbb Z\tag{2.14}
 \end{align*}
 where
 \begin{align*}
 \psi_{k,j,\lambda}^{M}(t)=(2N)^{j/2} \, \psi_k \big((2N)^j t - \lambda \big) \, e^{ - \frac{i\pi A}{B} ( t^2 -  \lambda^2 ) },~~ 1 \leq k \leq 2N-1, \lambda \in \Omega. \tag{2.15}
 \end{align*}
 Since the closed subspace $V_{1}^{M}$ can be decomposed as $V_{1}^{M}=V_{0}^{M}\oplus W_{0}^{M}$, so we have $\psi_{k}^{M}(t)\in W_{0}^{M}\subseteq V_{1}^{M}$, for $1\le k \le 2N-1$ and every fixed  $M=(A,B,C,D)$, and as a consequence, there exist a sequence $\big\{ b_{k,\lambda} \big \}_{ \lambda \in\Omega} $ with $\sum_{ \lambda \in \Omega} \big| b_{k,\lambda}\big|^2 < \infty $ such that
\begin{align*}
\psi_{k,0,0}^M(t)=\sqrt{2N} \sum_{\lambda\in\Omega} b_{k,\lambda} \,\phi\big(2N t-\lambda )e^{-\frac{i\pi A}{B}(t^2-\lambda^2 )},\tag{2.16}                            																		 \end{align*}
 which has an equivalent form in the LCT domain
 \begin{align*}
 \widehat{ \psi_k^M}\left( \frac{\omega}{B} \right) = \Lambda^M_k \left( \frac{\omega}{2NB} \right) \widehat{ \phi }\left( \frac{\omega}{2NB}\right),\tag{2.17}         													                             \end{align*}
where
\begin{align*}
 \Lambda^M_k\left( \frac{\omega}{B} \right) = \frac{1}{\sqrt{2N}} \sum_{ \lambda \in \Omega } b_{k,\lambda}^{M} \, e^{ -2 \pi i \lambda \omega/{B}}.\tag{2.18}                                                                               \end{align*}

\parindent=0mm \vspace{.0in}
Since $\Omega= \left\{0, r/N\right\}+2\mathbb Z$, therefore, we can write
\begin{align*}
  \Lambda^M_k\left( \frac{\omega}{B} \right) = \Lambda^{ M,1 }_k \left( \frac{\omega}{B} \right) + e^{ -2 \pi i \omega r/{NB}}\, \Lambda^{ M , 2 }_k \left( \frac{\omega}{B} \right), \tag{2.19}
  \end{align*}
where $\Lambda^{M,1}_k$ and $\Lambda^{M,2}_k$ are locally $ L^2$, $B/2$-periodic functions associated with the real uni-modular matrix $M=(A,B,C,D)$.

\parindent=8mm \vspace{.1in}
We are presently in a  position to build up the  completeness of the system $ \big \{\psi_{k,0,\lambda}(t)= \psi_k ( t - \lambda ) \, e^{ -i \pi{A}/{B} ( t^2 - \lambda^2 ) }, 1 \leq k \leq 2N-1 , \, \lambda \in \Omega\big \} $ in $ V^M_1 $ and subsequently,  we will able to obtain two orthonormality conditions for aforementioned system by means of $B/2$-periodic functions $ \Lambda^M_k$  given by (2.19).

\parindent=0mm \vspace{.1in}
{\bf Proposition 2.4.} {\it Consider a NUMRA associated with the real uni-modular matrix $M=(A,B,C,D)$ as in Definition 2.2. Suppose that there exist $(2N-1)$ functions $\psi^M_k, \, 1 \leq k \leq 2N - 1 $ in $ V^M_1 $. Then the family of functions}
\begin{align*}
\psi_{k,0,\lambda}^{M}(t)=\psi_k(x - \lambda ) \, e^{ - \frac{i \pi A}{B} ( t^2 - \lambda^2 ) } ,~~ 1 \leq k \leq 2N-1 , \, \lambda \in \Omega\tag{2.20}
\end{align*}
{\it will form an orthonormal system in $V^M_1$ if and only if}
\begin{align*}
 &\sum_{ p = 0 }^{ 2N - 1 } \left[ \Lambda^{M , 1}_\ell \left( \frac{\omega}{2NB} + \frac{p}{4N} \right) \overline{ \Lambda^{ M , 1 }_k \left( \frac{\omega}{2NB} + \frac{p}{4N} \right)} + \Lambda^{ M , 2 }_\ell \left( \frac{\omega}{2NB} + \frac{p}{4N} \right) \overline{ \Lambda^{ M , 2 }_k\left( \frac{\omega}{2NB} + \frac{p}{4N} \right) } \right] = \delta_{\ell,k }  \tag{2.21}  	
\end{align*}
{\it and}
\begin{align*}
   &\sum_{ p = 0}^{ 2N - 1 } e^{-i\pi r p/N} \left[ \Lambda^{M , 1}_\ell \left( \frac{\omega}{2NB} + \frac{p}{4N} \right) \overline{\Lambda^{M , 1}_k\left( \frac{\omega}{2NB} + \frac{p}{4N} \right)}  + \Lambda^{ M , 2 }_\ell \left( \frac{\omega}{2NB} + \frac{p}{4N} \right)  \overline{\Lambda^{ M , 2 }_k\left( \frac{\omega}{2NB} + \frac{p}{4N} \right) } \right] = 0\tag{2.22}
\end{align*}

\parindent=0mm \vspace{.0in}

{\it Proof.} Assume that the system given by (2.20) is orthonormal in $V^M_1$. Then, for $ \lambda, \sigma \in \Lambda$ and $0 \leq k ,\ell \leq 2N - 1 $, we have
\begin{align*}
  \Big \langle \psi_{\ell,0,\lambda}^{M}, \psi_{k,0,\sigma}^{M} \Big \rangle &=\Big \langle \psi_\ell ( t - \lambda ) \, e^{ -i \pi \frac{A}{B}( t^2 - \lambda^2 ) } , \psi_k ( t - \sigma ) \, e^{ -i \pi \frac{A}{B} ( t^2 - \sigma^2 ) } \Big \rangle   \\
 &= \int_{\mathbb R} \psi_\ell ( t - \lambda ) \, e^{ -i \pi \frac{A}{B} ( t^2 - \lambda^2 ) } \,\overline{ \psi_k( t - \sigma) } \, e^{ i \pi \frac{A}{B} ( t^2 - \sigma^2 ) } \, dt   \\
 &= e^{ i \pi \frac{A}{B} ( \lambda^2 - \sigma^2 ) } \int_{\mathbb R} \psi_\ell ( t - \lambda ) \,\overline{ \psi_k(t-\sigma) } \, dt   \\
 &= e^{ i \pi \frac{A}{B} ( \lambda^2 - \sigma^2 ) } \, \delta_{\ell,k } \, \delta_{ \lambda , \sigma }.
\end{align*}

In the LCT domain, we can write it as
\begin{align*}
 \delta_{\ell,k } \, \delta_{ \lambda , \sigma } = \frac{1}{B}\int_{\mathbb R}  \widehat{ \psi_\ell } \left(\frac{\omega}{B} \right) \,\overline{ \widehat{ \psi_k} \left( \frac{\omega}{B}\right) } \, e^{ - \frac{2 \pi i \omega}{B} ( \lambda - \sigma ) } \, d\omega.
\end{align*}
By taking $ \lambda = 2m , \, \sigma = 2n $, where $ m , n \in \mathbb Z $, we have
\begin{align*}
 \delta_{\ell,k } \, \delta_{ m , n } &= \frac{1}{B}\int_{\mathbb R}  \widehat{ \psi_\ell } \left( \frac{\omega}{B} \right) \, \overline{ \widehat{ \psi_k } \left( \frac{\omega}{B} \right) } \, e^{ -i4 \pi \frac{\omega}{B} ( m - n ) } \, d\omega  \\
 &= \frac{1}{B} \int_{ [ 0 ,B N ) } e^{-i4 \pi \frac{\omega}{B} ( m - n ) } \sum_{ j \in \mathbb Z } \widehat{ \psi_\ell } \left( \frac{\omega}{B} + N j \right) \, \overline{ \widehat{ \psi_k} \left( \frac{\omega}{B} + N j \right) } \,d\omega.
\end{align*}
Define
\begin{align*}
  h_{\ell,k } \left( \frac{\omega}{B} \right) = \sum_{ j \in \mathbb Z } \widehat{ \psi_\ell } \left( \frac{\omega}{B} + N j\right) \, \overline{ \widehat{ \psi_k} \left( \frac{\omega}{B} + N j \right) }.\tag{2.23}
\end{align*}
Then, we have
\begin{align*}
 \delta_{\ell,k }\,\delta_{ m , n } =\frac{1}{B}\int_{[0,BN)} e^{ - 4 \pi i {\omega}( m - n )/{B}} h_{\ell,k } \left(\frac{\omega}{B}\right)d\omega =\frac{1}{B} \int_{[0,\frac{B}{2})} e^{ - 4 \pi i {\omega}( m - n )/{B}}\sum_{p=0}^{2N-1} h_{\ell,k } \left( \frac{\omega}{B}  +   \frac{p}{2} \right)d\omega
\end{align*}
and
\begin{align*}
\quad \sum_{ p = 0 }^{ 2N - 1 } h_{\ell,k } \left( \frac{\omega}{B}  + \frac{p}{2} \right) = 2 \, \delta_{\ell,k }.\tag{2.24}
 \end{align*}
Again,  taking $ \lambda = r/N + 2m $ and $ \sigma = 2n $, where $ m , n  \in \mathbb Z $, we have
\begin{align*}
0 &= \frac{1}{B}\int_{\mathbb R}  e^{ -i2 \pi \frac{\omega}{B} \left( \frac{r}{N} + 2m - 2n \right) } \, \widehat{ \psi_\ell } \left( \frac{\omega}{B} \right) \, \widehat{ \psi_k } \left( \frac{\omega}{B} \right) d\omega  \\
&=  \frac{1}{B}\int_{ [ 0 , BN ) } e^{ -i4 \pi \frac{\omega}{B} ( m - n ) } \, e^{ -i2 \pi \frac{\omega}{B} \frac{r}{N} } \sum_{ j \in \mathbb Z } \widehat{ \psi_\ell} \left( \frac{\omega}{B} +  Nj \right) \, \overline{ \widehat{ \psi_k } \left( \frac{\omega}{B} + Nj \right) } \, d\omega      \\
&=  \frac{1}{B}\int_{ [ 0 ,B N ) } e^{ -i4 \pi \frac{\omega}{B} ( m - n ) } \, e^{ -i2 \pi \frac{\omega}{B} \frac{r}{N} }  \, h_{\ell,k } \left(  \frac{\omega}{B}\right) d\omega   \\
&=  \frac{1}{B}\int_{ \left[ 0 , \frac{B}{2} \right) } e^{ -4 \pi i \frac{\omega}{B} ( m - n ) } \, e^{ -i2 \pi \frac{\omega}{B} \frac{r}{N} }  \left[ \sum_{ p = 0 }^{ 2N - 1 } e^{ -i \pi \frac{r}{N} p } \, h_{\ell,k} \left( \frac{\omega}{B} + \frac{p}{2} \right) \right]d\omega .
\end{align*}

Therefore, we have
\begin{align*}
\sum_{ p = 0 }^{ 2N - 1 } e^{ -\pi ir p/{N} }\, h_{\ell,k } \left( \frac{\omega}{B} + \frac{p}{2} \right) = 0.\qquad\qquad\qquad\tag{2.25}
\end{align*}
The conditions given by (2.24) and (2.25) are equivalent to the orthonormality of the system given by (2.20). In terms of the high-pass filters $\Lambda_{k}^{M}$, conditions can be expressed as
\begin{align*}
 h_{\ell,k} \left( \frac{2N \omega}{B} \right) &= \sum_{ j \in \mathbb Z } \widehat{ \psi_\ell } \left( 2N \Big( \frac{\omega}{B} + \frac{j}{2} \Big) \right) \overline{ \widehat{ \psi_k} \left( 2N \Big( \frac{\omega}{B} + \frac{j}{2} \Big) \right) }  \\
 &= \sum_{ j \in \mathbb Z } \Lambda^M_\ell \left( \frac{\omega}{B} + \frac{j}{2}\right) \, \widehat{ \phi } \left( \frac{\omega}{B} +  \frac{j}{2}\right) \, \overline{ \Lambda^M_k \left( \frac{\omega}{B} + \frac{j}{2} \right) \, \widehat{ \phi } \left( \frac{\omega}{B} + \frac{j}{2} \right) }  \\
 &= \sum_{j \in \mathbb Z } \Lambda^M_\ell \left( \frac{\omega}{B} + \frac{j}{2}\right) \, \overline{\Lambda^M_\ell \left( \frac{\omega}{B} + \frac{j}{2} \right) } \, \left| \widehat{ \phi } \left( \frac{\omega}{B} + \frac{j}{2} \right) \right|^2  \\
 &= \left[\Lambda^{ M,1 }_\ell \left( \frac{\omega}{B} \right) \, \overline{\Lambda^{ M , 1 }_k \left( \frac{\omega}{B} \right) } + \Lambda^{ M , 2 }_\ell \left( \frac{\omega}{B} \right) \, \overline{ \Lambda^{M,2}_k \left( \frac{\omega}{B}\right) } \,\right] \sum_{ j \in \mathbb Z } \left| \widehat{ \phi } \left( \frac{\omega}{B} + \frac{j}{2}\right) \right|^2   \\
 & \quad + \left[ e^{ 2\pi i\frac{r\omega}{NB}} \Lambda^{M,1}_\ell \left( \frac{\omega}{B} \right) \, \overline{ \Lambda^{M,2}_k\left( \frac{\omega}{B} \right) } \sum_{ j \in \mathbb Z } e^{-ij\pi r/N} \left| \widehat{ \phi } \left( \frac{\omega}{B} + \frac{j}{2} \right) \right|^2 \right]    \\
  & \quad + \left[ e^{- 2\pi i\frac{r\omega}{NB} } \Lambda^{M,2}_\ell \left( \frac{\omega}{B} \right) \, \overline{ \Lambda^{M,1}_k\left( \frac{\omega}{B} \right) } \sum_{ j \in \mathbb Z } e^{ij\pi r/N} \left| \widehat{ \phi } \left( \frac{\omega}{B} + \frac{j}{2}\right) \right|^2 \right].
\end{align*}
Hence,
\begin{align*}
h_{\ell,k} \left(  \frac{2N\omega}{B} \right) &= \left[ \Lambda^{M,1}_\ell \left( \frac{\omega}{B} \right) \, \overline{ \Lambda^{M,1}_k \left( \frac{\omega}{B} \right) } + \Lambda^{M,2}_\ell \left( \frac{\omega}{B} \right) \, \overline{ \Lambda^{M,2}_k \left( \frac{\omega}{B} \right) }  \right] \sum_{j=0}^{ 2N - 1 } h_{ 0 , 0 } \left( \frac{\omega}{B} + \frac{j}{2} \right) \\
 & \quad + \left[ \Lambda^{M,1}_\ell\left( \frac{\omega}{B} \right) \, \overline{ \Lambda^{M,2}_k\left( \frac{\omega}{B} \right) } \, e^{2 \pi i \frac{r\omega}{NB} } \sum_{j=0}^{2N-1} e^{-ij\pi r/N} \, h_{0,0} \left( \frac{\omega}{B} + \frac{j}{2} \right) \right] \\
 & \quad + \left[ \Lambda^{M,2}_\ell \left( \frac{\omega}{B} \right) \, \overline{ \Lambda^{M,1}_k \left( \frac{\omega}{B} \right) } \, e^{ -2 \pi i \frac{r\omega}{NB}} \sum_{ j = 0 }^{ 2N - 1 } e^{ij\pi r/N} \, h_{ 0 , 0 } \left( \frac{\omega}{B} + \frac{j}{2} \right) \right]  \\
 &= 2 \left[ \Lambda^{M,1}_\ell \left( \frac{\omega}{B} \right) \, \overline{ \Lambda^{M,1}_k\left( \frac{\omega}{B} \right) } + \Lambda^{M,2}_\ell \left( \frac{\omega}{B} \right) \, \overline{ \Lambda^{M,2}_k \left( \frac{\omega}{B} \right) }  \right]. \tag{2.26}
\end{align*}
Applying the  identities (2.24)--(2.29), we obtain the desired conditions (2.21) and (2.22).

\parindent=0mm\vspace{.1in}
Now we shall  prove the sufficient part. From equations (2.17) and (2.26), we see that
\begin{align*}
&\sum_{ j \in \mathbb Z } \widehat{ \psi_\ell} \left( 2N \Big( \frac{\omega}{B}  +  \frac{j}{2} \Big) \right) \, \overline{ \widehat{ \psi_k} \left( 2N \Big( \frac{\omega}{B}  +  \frac{j}{2} \Big) \right) }  \qquad\qquad\qquad\\
&\qquad= \sum_{ j \in \mathbb Z } \Lambda^M_\ell \left( \frac{\omega}{B}  +  \frac{j}{2}\right) \, \overline{ \Lambda^M_k \left( \frac{\omega}{B}  +  \frac{j}{2}\right) } \, \left| \widehat{ \phi } \left( \frac{\omega}{B}  +  \frac{j}{2} \right) \right|^2  \\
&\qquad= 2 \left[ \Lambda^{M,1}_\ell \left( \frac{\omega}{2NB} \right) \, \overline{ \Lambda^{M,1}_k \left( \frac{\omega}{2NB}\right) } + \Lambda^{M,2}_\ell \left( \frac{\omega}{2NB} \right) \, \overline{ \Lambda^{M,2}_k \left( \frac{\omega}{2NB}\right) } \right]  \\
&\qquad= 2 \sum_{ p = 0 }^{ 2N - 1 } \left[ \Lambda^{M,1}_\ell \left( \frac{1}{2N} \left( \frac{\omega}{B} +  \frac{p}{2}\right) \right) \,  \overline{ \Lambda^{M,1}_k \left( \frac{1}{2N} \left( \frac{\omega}{B} + \frac{p}{2} \right) \right) } \right]   \\
&~~\qquad + \left[ \Lambda^{M,2}_\ell \left( \frac{1}{2N} \left( \frac{\omega}{B} +  \frac{p}{2}\right) \right) \, \overline{ \Lambda^{M,2}_k \left( \frac{1}{2N} \left( \frac{\omega}{B} +  \frac{p}{2} \right) \right) } \right] \\
&\qquad =  2 \, \delta_{\ell,k },
\end{align*}
which proves the orthonormality of the system (2.20). This completes the proof of Proposition 2.4. \qquad\fbox

\parindent=8mm \vspace{.0in}

Following result gives the existence of  nonuniform wavelet function associated with  LCT.

\parindent=0mm \vspace{.1in}
{\bf Theorem 2.5.} {\it Let $ \big \{ \psi_{k,0,\lambda}^{M} : 1 \leq k\leq 2N - 1 , \, \lambda \in \Omega \big \} $ be the system as defined in Proposition 2.4 and orthonormal in $ V^M_1 $. Then, this system is complete in $ W_{0}^{M}\equiv V^M_1 \ominus V^M_0$.}

\parindent=0mm \vspace{.1in}
{\it Proof.} The completeness of the system
\begin{align*}
\psi_{k,0,\lambda}^{M}(t)= \psi_k( t - \lambda ) \,  e^{ -i \pi \frac{A}{B} ( t^2 - \lambda^2 ) } ,~~1 \leq k \leq 2N - 1 , \, \lambda \in \Omega,\tag{2.27}
\end{align*}
is equivalent to the completeness of the system
\begin{align*}
\psi_{k,0,\lambda}^{M}(t)= \frac{1}{2N} \, \psi_k \left( \frac{t}{2N} - \lambda \right) e^{ -i \pi \frac{A}{B} ( t^2 - \lambda^2 ) },\quad 0 \leq k \leq 2N-1 , \, \lambda \in \Omega
\end{align*}
 in $ V^M_0 $. Therefore, under the given assumption, for every function $ f \in V^M_0 $, there exist a unique function $\Lambda^M_0 \big( {\omega}/{B} \big)$  of the form $ \frac{1}{\sqrt{ 2N }} \sum_{ \lambda \in \Omega } a^M_\lambda \, e^{- {2 \pi i\omega\lambda}/{B}} $ with $ \sum_{ \lambda \in\Omega} \big| a^M_\lambda \big|^2 < \infty $ such that
 \begin{align*}
  \widehat{f} \left(  \frac{\omega}{B}\right) =\Lambda^M_0 \left( \frac{\omega}{B}\right) \widehat{ \phi } \left( \frac{\omega}{B} \right).
  \end{align*}
 Therefore, it is enough to show that the family
\begin{align*}
\mathcal P = \Big \{ e^{ -2 \pi i (2N) \frac{\omega \lambda}{B}  } \Lambda^M_k\left( \frac{\omega}{B} \right) \chi_S \left( \frac{\omega}{B} \right) : 0 \leq k \leq 2N - 1 , \, \lambda \in \Omega\Big \}\tag{2.28}
\end{align*}
is complete in $ L^2 (S) $, where $ S\subset \mathbb R $ with $ 0 < |S| < \infty $.

\parindent=8mm \vspace{.1in}
Since the collection $ \big \{ e^{ -2 \pi i \frac{\omega \lambda }{B} } \, \chi_S \left( \frac{\omega}{B} \right) : \lambda \in \Omega\big \} $ is an orthonormal basis for $ L^2 (S ) $, therefore there exists locally $L^2$-functions $h_1 $ and $h_2 $ such that
\begin{align*}
h\left( \frac{\omega}{B} \right) = \Big[h_1 \Big( \frac{\omega}{B} \Big) + e^{ -2 \pi i \frac{r\omega}{NB} }h_2 \Big( \frac{\omega}{B} \Big) \Big] \chi_S \Big( \frac{\omega}{B} \Big).
\end{align*}

\parindent=0mm \vspace{.0in}
Assume that $h$ is orthogonal to all the functions belonging to the system (2.28),  then we see that
\begin{align*}
 0 &=\int_{S} e^{ -2 \pi i( 2N ) \frac{ \lambda \omega}{B}} \, \Lambda^M_k\left( \frac{\omega}{B} \right)  \overline{h \left( \frac{\omega}{B}\right) }\,d\omega \\
 &= \int_{ \left[ 0 , \frac{B}{2} \right) } e^{ -2 \pi i( 2N ) \frac{\lambda\omega}{B} } \left[ \Lambda^M_k\left( \frac{\omega}{B} \right) \, \overline{ h \left(\frac{\omega}{B}\right) } +  \Lambda^M_k \left( \frac{\omega}{B} + \frac{N}{2} \right) \overline{h\left( \frac{\omega}{B} + \frac{N}{2} \right) } \right] d\omega   \\
 &= \int_{ \left[ 0 , \frac{B}{2} \right) } e^{ -2 \pi i( 2N ) \frac{\lambda\omega}{B} } \left[ \Lambda^{M,1}_k \Big( \frac{\omega}{B} \Big) \, \overline { h_1 \left( \frac{\omega}{B} \right) } +   \Lambda^{M,2}_k \left( \frac{\omega}{B} \right) \, \overline{h_2 \left( \frac{\omega}{B}\right) } \right] d\omega .												 \end{align*}
For $\lambda = 2m$ and $k=0,1,\dots, 2N-1$, we define
\begin{align*}
 W_k \Big( \frac{\omega}{B} \Big) = \Lambda^{M,1}_k \left( \frac{\omega}{B}\right) \overline {h_1 \left( \frac{\omega}{B}\right) } +   \Lambda^{M,2}_k \left( \frac{\omega}{B} \right)  \overline{h_2 \left( \frac{\omega}{B} \right) },\tag{2.29}
\end{align*}
so that
\begin{align*}
0 &= \int_{ \left[ 0 , \frac{B}{2} \right) } e^{- 2 \pi i (4N)m \omega/{B} } \, W_k \left( \frac{\omega}{B} \right) d\omega\\
  &= \int_{ \left[ 0 , \frac{B}{4N} \right)} e^{- 2 \pi i (4N)m \omega/{B} } \sum_{ j = 0 }^{ 2N - 1 } W_k \left( \frac{\omega}{B} +  \frac{j}{4N} \right)d\omega .
\end{align*}

Since this equality holds for all $m\in\mathbb Z$, therefore
\begin{align*}
 \sum_{j=0}^{2N-1} W_k \left( \frac{\omega}{B} +  \frac{j}{4N} \right) = 0, \quad\text{a.e}. \tag{2.30}
\end{align*}
Similarly, on taking $\lambda =2m +{r}/{N}$, where $m\in\mathbb Z$, we have
\begin{align*}
 0 &= \int_{ \left[ 0 , \frac{B}{2} \right) } e^{ -i2 \pi (4N)m  \frac{\omega}{B}  } e^{ -i2 \pi (2r) \frac{\omega}{B} }  W_k \left( \frac{\omega}{B}\right)d\omega\\
 &= \int_{ \left[ 0 , \frac{B}{4N} \right) } e^{-i2 \pi (4N)m \frac{\omega}{B} } e^{-i4 \pi r \frac{\omega}{B}} \sum_{ j = 0 }^{ 2N - 1 } e^{-i\pi rj/N}  W_k \left( \frac{\omega}{B} + \frac{j}{4N} \right)d\omega.
\end{align*}

Hence, we deduce that
\begin{align*}
 \sum_{ j = 0 }^{ 2N - 1 }e^{-i\pi rj/N} \, W_k\left( \frac{\omega}{B} +  \frac{j}{4N} \right) = 0,\quad\text{a.e}.\tag{2.31}
\end{align*}

This completes the proof of Theorem 2.5. \qquad\fbox

\parindent=8mm\vspace{.1in}

If $ \psi^M_0 , \psi^M_1 ,\dots , \psi^M_{2N-1}  \in V^M_1 $ are as in Proposition 2.4, one can obtain from them  an orthonormal basis for $ L^2 ( \mathbb R ) $ by following the standard procedure for construction of wavelets from a given MRA \cite{M}. It can be easily checked that  for every $ j \in \mathbb Z $, the family $ \big \{ \psi_{k,j,\lambda}^{M}: 0\leq k \leq 2N-1, \lambda \in \Omega \big \} $ given by (2.15) constitutes a complete orthonormal system for $V_{j+1}^M$. Then, by virtue of (2.13),  it follows  that the system (2.15) constitutes an orthonormal basis for $L^2(\mathbb R )$.

\parindent=8mm\vspace{.1in}
In the following theorem, we establish a necessary and sufficient condition for the existence of basic wavelets associated with  LCT-NUMRA.

\parindent=0mm \vspace{.1in}

{\bf Theorem 2.6.} {\it Consider a  NUMRA  associated with the linear canonical transform as in Definition 2.2, such that the corresponding space $V_{0}^{M}$ has an orthonormal system of the form $\big\{\phi (t-\lambda ) \, e^{ -i \pi \frac{A}{B}(t^2-\lambda^2 ) } : \lambda \in \Omega\big\}$ and $\widehat{\phi}$ satisfies the two scale relation (2.8). Define}
\begin{align*}
\mathbb M_0 \left(\frac{\omega}{B}\right) = \left|\Lambda^{M,1}_0 \left( \frac{\omega}{B} \right)\right|^2 + \left| \Lambda^{M,2}_0 \left( \frac{\omega}{B}\right)\right|^2 ,\tag{2.32}
 \end{align*}
where $\Lambda^{M, 1 }_0 , \Lambda^{M , 2 }_0 $ are locally $L^2$-functions associated with the modular matrix  $M=(A,B,C,D)$. Then, a necessary and sufficient condition for the existence of associated linear canonical wavelets $\psi^M_1 ,\psi^M_2,\dots, \psi^M_{2N-1}$ is that $\mathbb M_0 $ satisfies the identity
\begin{align*}
 \mathbb M_0 \left( \frac{\omega}{B} + \frac{B}{4} \right) = \mathbb M_0 \left( \frac{\omega}{B}\right).\tag{2.33}
\end{align*}

\parindent=0mm \vspace{.0in}
{\it Proof.} The orthonormality of the system $\big\{\phi (t-\lambda ) \, e^{ -i \pi \frac{A}{B}(t^2-\lambda^2 ) } : \lambda \in \Omega\big\}$ which satisfies  (2.8), implies the following identities as shown in the proof of Proposition 2.5
\begin{align*}
&\sum_{ p = 0 }^{ 2N - 1 } \left[ \left|\Lambda^{M,1}_0 \left( \frac{\omega}{B} + \frac{p}{4N} \right)\right|^2 + \left|\Lambda^{M,2}_0 \left( \frac{\omega}{B} + \frac{p}{4N} \right) \right|^2 \right] = 1,\quad\text{and} \tag{2.34}\\
&\sum_{ p = 0 }^{ 2N - 1 } e^{-i\pi rp/N} \left[ \left| \Lambda^{M,1}_0 \left( \frac{\omega}{B} + \frac{p}{4N} \right) \right|^2 + \left| \Lambda^{M,2}_0 \left( \frac{\omega}{B} + \frac{p}{4N} \right) \right|^2 \right] = 0. \tag{2.35}
\end{align*}

Similarly, if $ \psi^M_k,k=0,1,\dots, 2N-1,$  are the basic  linear canonical nonuniform wavelets associated with the given NUMRA, then it satisfies the identity (2.17) and the orthonormality of the system  $ \big\{\psi^M_k: 0 \leq k \leq 2N - 1 \big \} $ in $ V^\alpha_1 $ is equivalent to the following identities:
\begin{align*}
\sum_{ p = 0 }^{ 2N - 1 } \left[ \Lambda^{M,1}_\ell \left( \frac{\omega}{B} +  \frac{p}{4N}  \right)  \overline{ \Lambda^{M,1}_k \left( \frac{\omega}{B} +  \frac{p}{4N}  \right) } +\Lambda^{M,2}_\ell \left( \frac{\omega}{B} +  \frac{p}{4N}  \right) \overline{ \Lambda^{M,2}_k \left( \frac{\omega}{B} +  \frac{p}{4N}  \right) }\right] = \delta_{\ell,k} \tag{2.36}
\end{align*}
and
\begin{align*}
 \sum_{ p = 0 }^{ 2N - 1 } e^{-i\pi rp/N} \left[ \Lambda^{M,1}_\ell \left( \frac{\omega}{B} +  \frac{p}{4N}  \right) \overline{ \Lambda^{M,1}_k\left( \frac{\omega}{B} +  \frac{p}{4N}  \right) }  +\Lambda^{M,2}_\ell \left( \frac{\omega}{B} +  \frac{p}{4N}  \right) \overline{ \Lambda^{M,2}_k \left( \frac{\omega}{B} +  \frac{p}{4N}  \right) } \right] = 0,\tag{2.37}
\end{align*}
where $  0 \leq k, \ell \leq 2N - 1 $. Moreover, if $ a_k( p ) = \Lambda^{M,1}_k \left( \frac{\omega}{B} +  \frac{p}{4N}  \right), \, b_k( p ) = \Lambda^{M,2}_\ell \left( \frac{\omega}{B} +  \frac{p}{4N}  \right)$  are vectors in $ \mathbb C^{ 2N } $ for $ p = 0,1,\dots ,2N - 1 , 0 \leq k \leq 2N - 1 $ , where $ u \in [ 0 , B/4 ] $ is fixed, then the solvability of the system of equations (2.36) and (2.37) is equivalent to
\begin{align*}
\mathbb M_0 \left( \frac{\omega}{B}+\frac{p+N}{4N}\right)=\mathbb M_0 \left( \frac{\omega}{B}+\frac{p}{4N}  \right), \quad \omega \in \left[ 0 , \frac{1}{4N}\right],
\end{align*}
$ p = 0, 1 , \dots , 2N - 1 $, which in turn is equivalent to  (2.33). The proof of this fact can be proved in similar lines as Lemma 3.5 in \cite{G1}. This completes the proof of Theorem 2.6. \qquad\fbox

\parindent=0mm\vspace{.1in}

{\bf Example 2.7.} For the demonstration of the above construction of nonuniform wavelets  associated with LCT, we  present an illustrative example of nonuniform Haar wavelets on the spectrum $\Omega= \{ 0, {1}/{N}\} + 2 \mathbb Z$.  For the parametric matrix $M=(A,B,C,D)$, we choose $\phi^M = \chi_{A_N}$, where
\begin{align*}
  A_N =  \bigcup_{j=0}^{N-1} \left[ \frac{2j}{N},\, \frac{2j + 1}{N} \right),\tag{2.38}
\end{align*}
is a measurable set of $\mathbb R$ with finite measure. Therefore, in view of specific form of $\Omega$, we can write
\begin{align*}
 \phi^M_{0,0} = \chi_{\left[ 0, {1}/{N} \right)}\circledast_M \sum_{j=0}^{N-1} \delta_{{2j}/{N}}.\tag{2.39}
\end{align*}
It is easy verify that the collection $\big\{ \phi_{0,\lambda}^{M}(t)= \phi(t - \lambda) \, e^{-i\pi \frac{A}{B} (t^2 - \lambda^2)} : \lambda\in\Lambda \big\}$ constitutes an orthonormal basis for $V^M_0$. Therefore, we have
\begin{align*}
 \frac{1}{2N} \phi^M_{0,0} \left( \frac{t}{2N}\right) &= \chi_{[0,2)}\circledast_M \sum_{j=0}^{N-1} \delta_{4j} = \left( \delta_0 + \delta_{\frac{1}{N}} \right)\circledast_M \phi^M_{0,0}\circledast_M \sum_{j=0}^{N-1} \delta_{4j},\tag{2.40}
  \end{align*}
which is required refinement equation for the given LCT-NUMRA. By taking the linear canonical transform on both sides of (2.40), we obtain
\begin{align*}
 \widehat{\phi} \left( \frac{2N\omega}{B} \right) = \frac{1}{2N} \left[ 1 + e^{i\pi \left( \frac{A}{BN^2} - \frac{2\omega}{BN} \right)} \right] \sum_{j=0}^{N-1} e^{i\pi \big( \frac{A(4j)^2}{B} \, - \, \frac{2(4j)\omega}{B} \big) }  \widehat{\phi} \left(\frac{\omega}{B}\right).\tag{2.41}
\end{align*}
Subsequently, we have
\begin{align*}
\Lambda^{M,1}_0 \left(\frac{\omega}{B}\right) = \Lambda^{M,2}_0 \left(\frac{\omega}{B}\right) = \frac{1}{2N}\sum_{k=0}^{N-1} e^{{i\pi A(4K)^2}/{B}}  e^{{-i \pi 8 \omega k}/{B}} .\tag{2.42}
\end{align*}
Clearly, condition (2.33) of the Theorem 2.6 is very well satisfied.

\parindent=0mm \vspace{.1in}
(i). For  $N = 1$ and any general  matrix $M=(A,B,C,D)$, we have $\phi^M_{0,0}(t)= \chi_{[0,1)}(t)$ and hence equation (2.42) yields $\Lambda^{M,1}_0 = \Lambda^{M,2}_0= 1/2$. Subsequently, the corresponding wavelet $\psi_1$ in the LCT domain is given by
\begin{align*}
  \widehat{\psi}_1\left(\frac{2\omega}{B}\right)= \frac{1}{2} \left( e^{-2\pi i{\omega}/{B}} - 1 \right)\widehat{\phi}\left(\frac{\omega}{B}\right).\tag{2.43}
\end{align*}
In particular, for the matrix $M=(2,1,1,1)$, we obtain linear canonical Haar wavelet given by
\begin{align*}
\psi_1(t)=
\begin{cases}
e^{-i\pi 8t^2}, \hspace{1.3cm} 0\leq t<{1}/{2}   \\
-e^{-i\pi 2 (2t-1)^2}, \quad {1}/{2} \leq t \leq 1,
   \end{cases}
\end{align*}
For the choice $M= (\cos\pi/4,\sin\pi/4,-\sin\pi/4,\cos\pi/4)$, we obtain the  fractional Haar wavelet of the form
\begin{align*}
  \psi_1(t)=
  \begin{cases}
e^{-i\pi 4t^2}, \hspace{1.3cm} 0\leq t<{1}{2}   \\
-e^{-i\pi 4 t(t-1)}, \quad {1}/{2} \leq t \leq 1.
\end{cases}
\end{align*}
Similarly, for the choice $M= (0,1,-1,0)$, one can obtain the classical Haar wavelet given by
\begin{align*}
  \psi_1(t)=
  \begin{cases}
1, \hspace{1.3cm} 0\leq t<{1}{2}   \\
-1, \quad {1}/{2} \leq t \leq 1.
\end{cases}
\end{align*}
(ii). For $N=2$ and parametric matrix $M=(A,B,C,D)$, the periodic linear canonical filters $\Lambda^{M,1}_0$ and $\Lambda^{M,2}_0$ given in (2.42) becomes
\begin{align*}
  \Lambda^{M,1}_0 \left(\frac{\omega}{B}\right) = \Lambda^{M,2}_0 \left( \frac{\omega}{B}\right)=\frac{e^{-i4\pi \big(\frac{2A}{B} + \frac{\omega}{B} \big)}}{2}   \cos \left( \frac{8\pi A}{B} +  \frac{4\pi\omega}{B} \right)
\end{align*}
so that the necessary and sufficient condition for the existence of associated wavelets given by relation (2.33) takes the form
\begin{align*}
\mathbb M_0\left(\frac{\omega}{B}\right)= \frac{1}{2} \cos^2 \left( \frac{8\pi A}{B} +  \frac{4\pi\omega}{B} \right).
\end{align*}
Consequently, the basic  wavelets  associated with given parametric matrix $M=(A,B,C,D)$ can be computed via relation (2.17) as
\begin{align*}
 &\widehat{\psi}_1 \left( \frac{4\omega}{B} \right) = \frac{i}{4} \left( 1+ e^{\frac{i\pi A}{B}} \right) e^{i \big( \frac{8\pi A}{B} - \frac{ 4\pi\omega}{B} \big) } \sin \left( \frac{8\pi A}{B} +  \frac{4\pi\omega}{B} \right) \left( 1+ e^{- \frac{i\pi\omega}{B}} \right),   \\
 &\widehat{\psi}_2 \left( \frac{4\omega}{B} \right) = \frac{1}{4} \left( 1+ e^{\frac{i\pi A}{B}} \right) e^{i \big( \frac{8\pi A}{B} + \frac{ 4\pi\omega}{B} \big) } \cos \left( \frac{8\pi A}{B} +  \frac{4\pi\omega}{B} \right) \left( -1+ e^{- \frac{i\pi\omega}{B}} \right),   \\
 &\widehat{\psi}_3 \left( \frac{4\omega}{B} \right) = \frac{i}{4} \left( 1+ e^{\frac{i\pi A}{B}} \right) e^{i \big( \frac{8\pi A}{B} + \frac{ 4\pi\omega}{B} \big) } \sin \left( \frac{8\pi A}{B} +  \frac{4\pi\omega}{B} \right) \left( -1+ e^{- \frac{i\pi\omega}{B}} \right).
\end{align*}
In particular, for the choice $M=(0,1,2,-1)$, we obtain the following set of basic wavelets in time domain
\begin{align*}
 &\psi_1(t) = \chi_{[0,1/2)} - \chi_{[1,3/2)}(t),    \\
 &\psi_2(t) = -\chi_{[-1,-7/8)}(t) + \chi_{[-7/8,-3/4)}(t) - \chi_{[-3/4,-5/8)}(t) + \chi_{[-5/8,-1/2)} (t) \\
  & \qquad \qquad - \chi_{[0,1/8)} (t)+ \chi_{[1/8,1/4)}(t) - \chi_{[1/4,3/8)}(t) + \chi_{[3/8,1/2)}(t),   \\
 & \psi_3(t) = -\chi_{[-1,-7/8)}(t) + \chi_{[-7/8,-3/4)}(t) - \chi_{[-3/4,-5/8)}(t) + \chi_{[-5/8,-1/2)} (t) \\
  & \qquad \qquad + \chi_{[0,1/8)}(t) - \chi_{[1/8,1/4)}(t) + \chi_{[1/4,3/8)}(t) - \chi_{[3/8,1/2)}(t).
\end{align*}

\section{An LCT-NUMRA Starting From A Scaling Function}

\parindent=0mm \vspace{.0in}
Recall that an LCT-NUMRA is a collection of nested spaces $\big\{V_{j}^{M}: j\in\mathbb Z\big\}$ of $L^2(\mathbb R)$ and a special function $\phi\in V_{0}^{M}$ such that (a)--(e) conditions of the Definition 2.2 are satisfied. However, it is also possible to construct an LCT-NUMRA by first choosing an appropriate scaling function $\phi(t)$ and obtaining $V_{0}^{M}$ by taking the linear span over the translations by elements of $\Omega$  of the given function $\phi$.  The other spaces $V_{j}^{M}$ can be generated as the scaled adaptations of $V_{0}^{M}$.

\parindent=8mm \vspace{.1in}
Let us  start with a function $\phi\in L^2(\mathbb R)$ such that
\begin{align*}
\Phi(t)=\phi (t) \, e^{ -i \pi \frac{A}{B} t^2 } = \sqrt{2N} \sum_{ \lambda \in \Omega} a_\lambda  \phi \big( 2Nt - \lambda \big) \, e^{- \frac{i \pi A}{B}( t^2 - \lambda^2 ) }, \tag{3.1}
\end{align*}
and $ \big\{\phi(t-\lambda)\,e^{-i\pi \frac{A}{B}(t^2-\lambda^2)}: \lambda\in\Omega\big\}$ is orthonormal in $L^2(\mathbb R)$. Then,  we define
\begin{align*}
V_{0}^{M}&=\overline{\text{Span}}\left\{\phi(t-\lambda)\,e^{-i\pi \frac{A}{B}(t^2-\lambda^2)}: \lambda\in\Omega\right\},\quad\text{and}\tag{3.2}\\
V_{j}^{M}&=\overline{\text{Span}}\left\{(2N)^{j/2} \phi\big((2N)^jt-\lambda\big)e^{-i\pi(t^2-\lambda^2)}: j \in \mathbb Z, \lambda\in\Omega\right\}.\tag{3.3}
\end{align*}
Condition (3.1) and orthonormality of the system $\big\{\phi_{0,\lambda}^{M}: \lambda\in \Omega\big\}$  guarantee that $\big\{\phi_{j,\lambda}^{M}: j\in\mathbb Z,\lambda\in \Omega\big\}$ is an orthonormal basis in each $V_{j}^{M}$ and hence, $V_{j}^{M}$ satisfy the increasing property $V_{j}^{M}\subset V_{j+1}^{M}$ of an LCT-NUMRA. Besides, it follows that $V_{j}^{M}$ satisfy scaling and translating properties of Definition 2.2 also. Now, in order to verify that the ladder of spaces generated by $\phi$ constitutes an LCT-NUMRA, it is sufficient to show that the following properties also hold:
\begin{align*}
\bigcap_{j\in\mathbb Z} V_{j}^{M}=\left\{ 0\right\}, \quad\text{and}\quad \overline{\bigcup_{j\in\mathbb Z} V_{j}^{M}}=L^2(\mathbb R).
\end{align*}
We shall first show that if $\phi(t)$ is of the form (2.12) and satisfies the following conditions
\begin{align*}
 \sum_{ p = 0  }^{ 2 N - 1 } \mathbb M_0 \left( \frac{\omega}{B} + \frac{p}{4N} \right) = 1,\quad\text{and}\quad \sum_{ p = 0 }^{ 2N - 1 } e^{-i\pi rp/N} \, \mathbb M_0 \left( \frac{\omega}{B} + \frac{p}{4N} \right) = 0,\tag{3.4}
\end{align*}
where $\mathbb M_{0}(\omega)$ is given by (2.32), then $\phi$ belongs to $L^2(\mathbb R)$.

\parindent=0mm \vspace{.1in}

{\bf Theorem 3.1.} {\it Let $\Lambda^M_0 (\omega)$ be of the form (2.8) and $\mathbb M_0(\omega)$ satisfying (3.4). Let $\phi(t)$ be defined by (2.12) and assume that the infinite product defining $ \widehat{\phi}(\omega) $ converges a.e on $\mathbb R$. Then, the function $\phi(t)$ belongs to $L^2(\mathbb R)$.}

\parindent=0mm \vspace{.1in}
{\it Proof.} Consider the integrals
\begin{align*}
  R_1 &= \int_{ |\omega| \leq \frac{N}{2} } \mathbb M_0 \left( \frac{\omega}{2NB} \right)d\omega,\\
  R_J &= \int_{ |\omega| \leq \frac{N(2N)^{J-1}}{2}} \mathbb M_0 \left( \frac{\omega}{(2N)^J B} \right)\prod_{k=1}^{J-1}\left|\Lambda^M_0 \left( \frac{\omega}{(2N)^k B}  \right) \right|^2 d\omega, \quad J \geq 2.
\end{align*}
Then, we have
\begin{align*}
  R_1 &= \int_{ [ 0 , N ) } \mathbb M_0 \left( \frac{\omega}{(2N)^{J}B} \right)d\omega= \int_{\left[0,\frac{B}{2}\right)} \sum_{p=0}^{2N-1}  \mathbb M_0 \left( \frac{\omega}{2NB} + \frac{p}{4N} \right)d\omega= \dfrac{1}{2},\\
  R_{J+1} &= \int_{ [0,N(2N)^J]} \mathbb M_0 \left( \frac{\omega}{(2N)^{J+1}B}\right) \prod_{k=1}^J \left|\Lambda^M_0 \left(\frac{\omega}{(2N)^k B}\right)\right|^2 d\omega\\
  &= \int_{ [0 ,N(2N)^{J-1}]} \sum_{p=0}^{2N-1} \left| \Lambda^M_0 \left( \frac{\omega}{B} +  \frac{p}{2}  \right) \right|^2 \mathbb M_0 \left( \frac{\omega}{2NB} +  \frac{p}{4N} \right) \prod_{k=1}^{J-1} \left| \Lambda^M_0 \left( \frac{\omega}{(2N)^k B} \right) \right|^2d\omega.
\end{align*}

\parindent=0mm\vspace{.0in}
We observe that
 \begin{align*}
 &\sum_{ p = 0 }^{ 2N - 1 } \left| \Lambda^M_0 \Big( \frac{\omega}{B} + \frac{p}{2} \Big) \right|^2 \mathbb M_0 \left( \frac{\omega}{2NB} \right)\\
 &\qquad= \left[ \left| \Lambda^{M,1}_0 \left( \frac{\omega}{B}\right) \right|^2 + \left| \Lambda^{M,2}_0 \left( \frac{\omega}{B} \right) \right|^2 \right] \sum_{ p = 0 }^{ 2N - 1 } \mathbb M_0 \left( \frac{\omega}{2NB} + \frac{p}{4N} \right) \\
 &~~\qquad+ \left[ \overline{ \Lambda^{M,1}_0 \left( \frac{\omega}{B} \right) }  \Lambda^{M,2}_0 \left( \frac{\omega}{B} \right) \, e^{ -i2 \pi \frac{\omega}{B}  \frac{r}{N}} \sum_{ p = 0 }^{ 2N - 1 } e^{-i\pi rp/N} \, \mathbb M_0 \left( \frac{\omega}{2NB} + \frac{p}{4N} \right)   \right] \\
 &~~\qquad+ \left[ \Lambda^{M,1}_0 \left( \frac{\omega}{B} \right) \, \overline{ \Lambda^{M,2}_0 \left( \frac{\omega}{B} \right) } \, e^{ i2 \pi \frac{\omega}{B}  \frac{r}{N}} \sum_{ p = 0 }^{ 2N - 1 } e^{i\pi rp/N} \mathbb M_0 \left( \frac{\omega}{2NB} + \frac{p}{4N} \right) \right]     \\
 &\qquad= \left| \Lambda^{M,1}_0 \left( \frac{\omega}{B} \right) \right|^2 + \left| \Lambda^{M,2}_0 \left( \frac{\omega}{B} \right) \right|^2 \\
&\qquad=  \mathbb M_0 \left( \frac{\omega}{B} \right).
\end{align*}
This shows that
\begin{align*}
R_{ J + 1 } = \int_{ [ 0 , N (2N)^{ J - 1 } ] } \mathbb M_0 \left( \frac{\omega}{(2N)^J B} \right) \prod_{ k = 1 }^{ J - 1 } \left| \Lambda^M_0 \left( \frac{\omega}{(2N)^k B} \right) \right|^2 d\omega= R_J, ~ \text{ for any } \, J \geq 1 .
\end{align*}
Therefore, it follows by induction  that
\begin{align*}
R_{ J + 1 } =R_J = R_{ J - 1 } = \dots = R_2 = R_1 = \dfrac{1}{2}.
\end{align*}
Hence,
\begin{align*}
\int_{ |\omega| \leq \frac{ N(2N)^{J - 1}}{2} } \left| \widehat{ \phi } \left( \frac{\omega}{B}\right)\right|^2 d\omega  & \leq \int_{ |\omega| \leq \frac{N(2N)^{J - 1}}{2} } \left| \Lambda^M_0 \left( \frac{\omega}{(2N)^J B} \right) \right|^2 \prod_{ k = 1 }^{ J -  1 } \left| \Lambda^M_0 \left( \frac{\omega}{(2N)^k B} \right) \right|^2d\omega\\
 & \leq \int_{ |\omega| \leq \frac{N (2N)^{J - 1}}{2} } 2 \mathbb M_0 \left( \frac{\omega}{(2N)^J B} \right) \prod_{ k = 1 }^{ J - 1 } \left|\Lambda^M_0 \left( \frac{\omega}{(2N)^k B} \right)  \right|^2d\omega\\
 & = 2 R_{ J - 1 }  \\
 &= 1.
\end{align*}
Since $J$ is arbitrary, it follows that $\phi\in L^2(\mathbb R)$. This completes the  proof of Theorem 3.1.  \qquad\fbox

\parindent=8mm\vspace{.1in}

For each $j\in\mathbb Z$ and fixed unimodular matrix $M=(A,B,C,D)$, we define the orthogonal projection $P_{j}^{M}$  of $L^2(\mathbb R)$ onto $V_{j}^{M}$ as
\begin{align*}
 P_{j}^{M} f = \sum_{\lambda \in\Omega} \left \langle f , \phi^M_{j,\lambda } \right \rangle \phi^M_{j,\lambda}. \tag{3.5}
\end{align*}

\parindent=0mm \vspace{.0in}

{\bf Theorem 3.2.} {\it Let $\big\{ V_{j}^{M}: j\in\mathbb Z\big\}$ be the family of subspaces as defined in (3.3) with given $\phi\in L^2(\mathbb R)$. If $\big \{ \phi(t-\lambda ) e^{-i\pi\frac{A}{B}(t^2-\lambda^2)} :\lambda \in\Omega\big\} $ is an orthonormal basis in $ V^M_0$, then $\bigcap_{j\in\mathbb Z } V^M_j=\left\{0 \right\}$.}

\parindent=0mm \vspace{.1in}
{\it Proof.}  Let $f\in\bigcap_{j\in \mathbb Z } V^M_j$. Then for $ \epsilon > 0 $, there exist compactly supported continuous function $f_\varepsilon $ such that $\big\| f - f_\varepsilon\big\| < \varepsilon $. Since $P_{j}^{M}$ is the orthogonal projection operator given by the formula (3.5), then we have
\begin{align*}
\left\|f-P^M_j f_\varepsilon\right\|_2 = \left\|P^M_j(f-f_\varepsilon) \right\|_2 \leq \big\|f-f_\varepsilon\big\|_2 < \varepsilon
\end{align*}
which implies that $\big\| f \big\|_2 < \varepsilon + \big\|P^M_j f_\varepsilon\big\|_2.$ Using the fact that the collection $\big\{ \phi_{j,\lambda}^{M}: j\in\mathbb Z, \lambda\in \Omega\big\}$ forms an orthonormal basis for $ V^M_j$, therefore, for each $j\in\mathbb Z$, we have
\begin{align*}
\left\|P^M_j f_\varepsilon\right\|^2_2 &= \sum_{ \lambda \in \Lambda } \left| \big \langle P^M_j f_\varepsilon, \phi^M_{j,\lambda}\big \rangle\right|^2 \\
&=\sum_{ \lambda \in \Omega } \left| \left \langle f_\varepsilon,\phi^M_{j,\lambda} \right \rangle \right|^2\\
&= (2N)^j \sum_{ \lambda \in \Omega } \left| \int_{-1/4}^{1/4} f_\varepsilon( t ) \phi \left((2N)^j t -\lambda \right) \, e^{i\pi\frac{A}{B}(t^2 - \lambda^2 ) }  dt \right|^2  \\
&\leq \frac{(2N)^j}{2} \sum_{ \lambda \in \Omega } \Big[ \big| f_\varepsilon ( t ) \big| \big| \phi \left( (2N)^j t - \lambda \right) \big| \, dt \Big]^2 \\
&\leq \frac{(2N)^j K^2}{2} \sum_{ \lambda \in \Omega } \Big[ \big| \phi \left( (2N)^j t - \lambda\right) \big|  dt \Big]^2 \\
&= \frac{K^2}{2}  \int_{\cup_{ \lambda \in \Omega } \left[ - \lambda- \frac{(2N)^j}{4}  ,  - \lambda + \frac{(2N)^j}{4} \right]} \chi_{ \cup_{ \lambda \in \Lambda } \left[ - \lambda- \frac{(2N)^j}{4} ,- \lambda + \frac{(2N)^j}{4} \right]}  \big| \phi(y) \big|^2 dy,
\end{align*}
where $K=\sup_{x\in\mathbb R}| f_\varepsilon (t)|$. By Dominated Convergence Theorem, this term tends to zero as $j \rightarrow \infty$. In particular, there exist a $j$ such that $\left\| P^M_j f_\varepsilon\right\|\rightarrow 0$ and therefore,  $\big\|f\big\|_2<\varepsilon$. Since $\varepsilon > 0$ is arbitrary, this implies that $f=0$, and hence, $\bigcap_{j\in \mathbb Z}V^\alpha_j =\left\{0\right\}$. This completes the proof of Theorem 3.2. \qquad\fbox

\parindent=0mm \vspace{.1in}

{\bf Theorem 3.3.} { \it  Let $\phi\in L^2( \mathbb R)$ such that $\big\{ \phi(t -\lambda) e^{-i\pi\frac{A}{B}(t^2-\lambda^2)}: \lambda \in\Omega \big\} $ is an orthonormal basis in $ V^M_0 $ and let $ \big \{ V^M_j:j\in\mathbb Z \big \}$ be the family of subspaces defined by (3.3). Assume that $\widehat {\phi}(\omega)$ is bounded and continuous near $\omega=0$ with $\big | \widehat { \phi }( 0 ) \big | \neq 0 $. Then, $ \overline { \bigcup_{j \in \mathbb Z } V^M_j}=L^2 ( \mathbb R )$. }

\parindent=0mm \vspace{.1in}
{\it Proof.} Let $ f \in \big ( \bigcup_{j \in \mathbb Z } V^M_j \big )^\perp $. Then, for $ \varepsilon > 0 $, there exists a $C^{\infty}$-function $ f_\varepsilon $ with compact support such that $ \big\| f - f_\varepsilon \big\| < \varepsilon $. Therefore
\begin{align*}
  &\big \| P_j^M f \big \|^2_2 = \big \langle P_j^M f , P_j^M f \big \rangle_2 = \big \langle f , P_j^M f \big \rangle_2 = 0,\quad\text{and}\\
  &\big \| P_j^Mf_\varepsilon \big \|^2_2 = \big \| P_j^M ( f - f_\varepsilon ) \big \|_2 \leq  \big \| f - f_\epsilon \big \|_2 < \varepsilon .
\end{align*}
Since the system $ \big \{ \phi_{j,\lambda}^{M}:j\in\mathbb Z, \lambda\in\Lambda\big\}$ constitutes an orthonormal basis for $ V^M_j $ and $ f_\varepsilon $ is of compact support, therefore we have
\begin{align*}
\left\| P_j^M f_\varepsilon \right\|^2_2 & = \sum_{ \lambda \in \Omega }\left| \big \langle f_\varepsilon , \phi^M_{j , \lambda}\big \rangle\right|^2 = \sum_{\lambda \in \Lambda} \left | \int_{\mathbb R} \mathscr L_M [f_\varepsilon](\omega)  \overline { \mathcal K_M \left( \frac{\lambda }{( 2N )^j} , \omega \right ) \widehat { \phi } \left (  \frac{\omega}{(2N)^j B} \right ) } d\omega \right |^2.\tag{3.6}                                                                                \end{align*}
Next, we choose $ j $ sufficiently large so that $ \text { supp } f_\varepsilon \subseteq  [ -1/4 , 1/4  ] $ and for this choice of $ j $, we assume
\begin{align*}
 H\Big( \frac{\omega}{B} \Big) = \mathscr L _M\big[f_\varepsilon\big](\omega) \, \overline { \widehat { \phi } \left( \frac{\omega}{( 2N )^j B} \right) },\tag{3.7}
\end{align*}
for some function $H$ of the form
\begin{align*}
H\left( \frac{\omega}{B}\right) = H^M_1 \left( \frac{\omega}{B}\right) + e^{2 \pi i\frac{\omega}{B} \frac{r}{N} }H^M_2 \left( \frac{\omega}{B}\right),\tag{3.8}
\end{align*}
where $ H^M_1 $ and $H^M_2 $ are locally square integrable periodic functions associated with the unimodular matrix $M=(A,B,C,D)$. Moreover, if $H(\omega)$ has the expansion of the form $\sum_{ \lambda \in \Omega} b^M_\lambda e^{-2\pi i\lambda {\omega}/{B} }$ on the set $S=  [ 0 , 1/2  ) \cup [N/2 , (N+1) / 2 ) $, then
\begin{align*}
  b^M_\lambda = \int_{S} H \Big( \frac{\omega}{B} \Big) \, e^{-2 \pi i \frac{\omega}{B} \lambda } d\omega= \int_{\mathbb R} \mathscr L_M \big[f_\varepsilon\big] \big( (2N)^j\omega\big) \overline {\widehat {\phi}\left(\frac{\omega}{B}\right) } e^{2 \pi i\frac{\omega}{B} \lambda } d\omega , \quad \lambda \in \Omega .
\end{align*}
Taking $\lambda=2m$, where $m\in\mathbb Z$, we have
\begin{align*}
 \int_{ \left[ 0 , \frac{B}{2} \right) } 2 \, H^M_1 \left(\frac{\omega}{B}\right) \, e^{2 \pi i(2k) \frac{\omega}{B} }\, d\omega &= \int_{ \left[ 0 , \frac{B}{2} \right) } \sum_{k\in \mathbb Z} \mathscr L_M \big[f_\epsilon\big] \left ((2N)^j\omega + \frac{k (2N)^j}{2} \right )\overline{\widehat{\phi}\left(\frac{\omega}{B} + \frac{k}{2} \right) } \, e^{2\pi i(2k) \frac{\omega}{B} }\, d\omega .
\end{align*}
Therefore,
\begin{align*}
 H^M_1 \Big( \frac{\omega}{B} \Big) = \dfrac{1}{2} \sum_{k\in \mathbb Z} \mathscr L_M \big[f_\epsilon\big] \left ((2N)^j\omega +\frac{k(2N)^j}{2}\right)\overline { \widehat{\phi}\left(\frac{\omega}{B} + \frac{k}{2} \right )}.
\end{align*}
Similarly, on taking $ \lambda = 2m + r/N$, where $m\in\mathbb Z$, we obtain
\begin{align*}
H^M_2 \left(\frac{\omega}{B}\right)=\dfrac{1}{2}\sum_{k\in\mathbb Z }\mathscr L_M \big[f_\epsilon\big] \left ((2N)^j\omega+\frac{k(2N)^j} {2} \right) \overline { \widehat {\phi}\left (\frac{\omega}{B}+\frac{k}{2} \right)} \, e^{2\pi i \left ( \frac{\omega}{B} + \frac{k}{2}\right)\frac{r}{N}} .
\end{align*}
Therefore, it follows that
\begin{align*}
H\left(\frac{\omega}{B}\right) =  \dfrac{1}{2} \sum_{k\in\mathbb Z} \mathscr L_M \big[f_\epsilon\big] \left((2N)^j \omega+ \frac{k(2N)^j} {2}\right) \overline { \widehat {\phi} \left(\frac{\omega}{B}+\frac{k}{2}\right)} \left (1+e^{-i\pi rk/N} \right ) .
\end{align*}
Since $ \text{supp} f_\varepsilon \subseteq [-1/4,1/4 ]$, therefore for large values of $j$, (3.6) becomes
\begin{align*}
 \big \| P_j^M f_\epsilon \big \|^2_2 \leq \int_{ \cup_{ j \in \mathbb Z } \left[-\frac{1}{4} + Nj ,\frac{1}{4}+ Nj \right]} \left|\mathscr L_M \big[f_\epsilon\big](B\eta)\overline{ \widehat{\phi} \left (\frac{\eta}{(2N)^j} \right)}\right|^2 d\eta .
\end{align*}

By Dominated Convergence Theorem, we observe that the R.H.S of above inequality converges to $ |\widehat{\phi}(0)|^2 \| f_\varepsilon \|^2_2 $, as $j\to\infty$. Therefore, we have
\begin{align*}
 \varepsilon > \big \|P^M_j f_\varepsilon \big \|^2_2  = \big \| \mathscr L_M [f_\varepsilon] \big \|^2_2 = \big \|f_\varepsilon \big \|^2_2.
\end{align*}
Consequently, $ \big \| f \big \|_2 < \varepsilon + \big \| f_\varepsilon \big \|_2 < 2 \varepsilon $. Since $\varepsilon$ is arbitrary, therefore $ f = 0 $. This completes the proof of Theorem 3.3. \qquad\fbox

\section{Nonuniform Wavelet Packets Associated with LCT}

\parindent=0mm \vspace{.0in}
In this section, we construct nonuniform wavelet packets associated with linear canonical transform and investigate their  orthonormal properties.

 \parindent=8mm\vspace{.1in}

Let $\big\{V_{j}^M:j\in\mathbb Z\big\}$ be an NUMRA associated with the linear canonical transform with scaling function $\phi$ satisfying (2.8), where $\Lambda_{0}^{M}$ is given by (2.10). Contrary to the standard procedure for the construction of wavelets, the existence of the  functions $\Lambda_{k}^{M}, 1\le k\le 2N-1$ satisfying (2.21) and (2.22) is not guaranteed. However, we proved in Section 2 that a necessary and sufficient condition for the existence of these locally $L^2$-functions (basic wavelets), so that (2.21) and (2.22) are satisfied is that the function $\mathbb M_{0}$, defined by (2.32) is $B/4$-periodic. As we know that the general technique involved for the construction of wavelet packets is through splitting the wavelet spaces $W_{j}^M$ successively into a finite number of orthogonal sub-spaces \cite{CM}. Since  the collection
\begin{align*}
\psi_{k,\lambda}^{M}(t)= \psi_k(t - \lambda)\, e^{-i\pi \frac{A}{B}(t^2-\lambda^2)},\quad 0 \leq k\leq 2N-1,\lambda\in\Omega,
 \end{align*}
is orthonormal if and only if $\Lambda^{M,1}_k$ and $\Lambda^{M,2}_k, 1\le k\le 2N-1$, defined by (2.10) satisfy (2.21) and (2.22) (Proposition 2.4). In deed, this
system is complete in $ V^M_{1}$ if and only if it is orthonormal (Theorem 2.5).

\parindent=8mm\vspace{.1in}
Applying this  technique to the space $V^M_1$, we obtain functions $\mathcal W_k, 0 \leq k\leq 2N-1 $, where
\begin{align*}
 \widehat {\mathcal W}_k\left(\frac{\omega}{B}\right)=\Lambda^M_{k} \left(\frac{\omega}{2NB}\right) \widehat{\phi}\left(\frac{\omega}{2NB}\right),																										 \tag{4.1}
\end{align*}

\parindent=0mm \vspace{.0in}
such that $\big \{\mathcal W_k(t-\lambda) \, e^{-i \pi \frac{A}{B} (t^2-\lambda^2)} : 0\leq k \leq 2N-1, \lambda\in\Omega\big\}$ forms an orthonormal basis for $V^M_1$. For $k=0$, we obtain the scaling function $\mathcal W_0= \phi$, whereas for  $1\le k\le 2N-1$, we have the basic wavelets $\mathcal W_{k}= \psi_k^{M},\,k= 1,2, \dots, 2N-1$. We now define $\mathcal W_{n}$ for each integer $n\ge 0$ as follows.

\parindent=0mm \vspace{.1in}
{\bf Definition 4.1.} The family of  nonuniform functions $\big\{ \mathcal W_{2Nn+k}(t):n\ge 0,0\le k\le 2N-1 \big\}$ is called a  nonuniform wavelet packet associated with a parameter matrix $M=(A,B,C,D)$ and  orthogonal scaling function $\phi$, where
\begin{align*}
\mathcal W_{2Nn+k}(t)=\sqrt{2N}\sum_{\lambda\in\Omega} b_{k,\lambda}^{M} \mathcal W_n (2Nt - \lambda) \, e^{-i \pi \frac{A}{B} (t^2-\lambda^2)},\quad 0 \leq k \leq 2N-1. \tag{4.2}
\end{align*}
Taking LCT on both sides of (4.2), we obtain
\begin{align*}
\widehat{\mathcal W}_{2Nn+\ell}\left(\frac{\omega}{B}\right)=\Lambda^M_k\left(\frac{\omega}{2NB}\right)\widehat{\mathcal W}_n \left(\frac{\omega}{2NB}\right), \quad   0 \leq k \leq 2N-1.\tag{4.3}
\end{align*}

The functions $\big\{ \mathcal W_{n}: n\ge 0\big\}$ will be called the basic wavelet packets associated with the given  linear canonical NUMRA or, simply, basic   nonuniform canonical wavelet packets.

\parindent=0mm\vspace{.1in}

Definition 4.1 allows us to make the following comments:

\begin{itemize}

  \item For the matrix $M=\left(0,1,-1,0\right)$, the Definition 4.1 boils down to the ordinary nonuniform wavelet packets \cite{Beh}.

  \item The classical wavelet packets \cite{CM} are obtained by obtained by switching the usual matrix $M=(A,B,C,D)$ to $M=\left(0,1,-1,0\right)$ and $N=1$.

  \item For $N=1$ and  $M=\left(\cos\theta,\sin\theta,-\sin\theta,\cos\theta\right)$, $\theta\neq n\pi$, we can obtain novel fractional wavelet packets as given by Shah and Debnath \cite{Shah6}.

  \item Similarly, for $M=(1,B,0,1),B\neq 0$, we can obtain a new class of wavelet packets called the nonuniform Fresnel-wavelet packets, whereas for the choice $M=\left(1,iB,0,1\right)$, we may have the  nonuniform Gauss-Weierstrass wavelet packets.

\end{itemize}

\parindent=0mm\vspace{.0in}
Moreover, if
\begin{align*}
n=\mu_{1}+(2N)\mu_2+(2N)^2\mu_{3}+\dots+(2N)^{j-1}\mu_{j},\tag{4.4}
\end{align*}
where $ 0\le\mu_i\le2N-1,1\le i\le j, \mu_{j}\ne 0,$ is the unique expansion of the non-negative integer $n$ in base $2N$. Then, an easy computation shows that the basic wavelet  packets ${\mathcal W_{n}}$ can expressed  in terms of the periodic  functions $\Lambda_{\mu_j}^M$ in the LCT domain as
\begin{align*}
\widehat{\mathcal W}_n\left(\frac{\omega}{B}\right)=\prod_{j=0}^{\infty}\Lambda_{\mu_j}^M\left(\frac{\omega}{(2N)^{j+1}B}\right).\tag{4.5}
\end{align*}

\parindent=0mm\vspace{.0in}
Further, it is worth noticing that the equations (2.24) and (2.25) can also be written as
\begin{align*}
&\sum_{p=0}^{2N-1} h_n \left(\frac{\omega}{B} + \frac{p}{2}\right) = 2, \quad \text{and}\quad \sum_{p=0}^{2N-1} e^{-\pi irp/N }\, h_n\left(\frac{\omega}{B} + \frac{p}{2}\right)=0,\tag{4.6}
\end{align*}
where
\begin{align*}
 h_n \left(\frac{\omega}{B}\right) = \sum_{j\in\mathbb Z} \widehat{\mathcal W}_n \left(\frac{\omega}{B} + Nj \right) \overline{\widehat{\mathcal W}_n \left( \frac{\omega}{B} + Nj \right)}.\qquad\qquad\qquad\tag{4.7}
\end{align*}

\parindent=0mm\vspace{.0in}
We are now ready to investigate the properties of the LCT-nonuniform wavelet packets by means of linear canonical transforms. First, we shall show that the proposed LCT nonuniform wavelet packets preserve the shift orthogonality property of the LCT  scaling function $\mathcal W_{0}(t)=\phi(t)$.

\parindent=0mm\vspace{.1in}
{\bf Theorem 4.2.} {\it Let $\phi(t)$ be a given LCT scaling function and $\big\{\mathcal W_{n}:n\ge 0\big\}$ be its corresponding family of nonuniform wavelet packets. Then, for each $n\in\mathbb Z^+$, we have}
\begin{align*}
  \Big\langle \mathcal W_n (t- \lambda),\mathcal W_n (t- \sigma) \Big\rangle = \delta_{\lambda,\sigma}, \quad \lambda,\sigma \in \Omega. \tag{4.8}
\end{align*}

\parindent=0mm \vspace{.0in}
{\it Proof.} We shall prove this result by induction on $n$. If $n=0$, then the result follows directly from the hypothesis. Assume that (4.8) holds for $ 0\leq n \leq (2N)^k$, where $k$ is a fixed positive integer. For $(2N)^k\leq n \leq (2N)^{k+1}$, we have  $(2N)^{k- 1} \leq \left[ \frac{n}{2N} \right] \leq (2N)^k$, where $[t]$ is the greatest integer function of $t$ and order $n=2N\left[ \frac{n}{2N} \right]+\ell,\, \ell=0,1,2,\dots, 2N-1$. By induction, we have
\begin{align*}
\left\langle \mathcal W_{\left[ \frac{n}{2N} \right]} (t- \lambda), \mathcal W_{\left[ \frac{n}{2N} \right]} (t- \sigma) \right\rangle =\delta_{\lambda,\sigma}, \quad \lambda,\sigma \in \Omega. \tag{4.9}
\end{align*}
By virtue of equation (4.3), we have
\begin{align*}
 &\Big\langle \mathcal W_n (t- \lambda),\mathcal W_n (t- \sigma) \Big\rangle \\
 & = \frac{e^{-i \pi \frac{A}{B}(\lambda^2 - \sigma^2)}}{B} \int_{\mathbb R} e^{-2\pi i\frac{\omega}{B} (\lambda - \sigma)} \,\widehat{\mathcal W}_{\left[ \frac{n}{2N} \right]}\left(\frac{\omega}{B}\right) \overline{\widehat{\mathcal W}_{\left[ \frac{n}{2N} \right]}\left(\frac{\omega}{B}\right)}\,d\omega\\
 & = \frac{e^{-i \pi \frac{A}{B}(\lambda^2 - \sigma^2)}}{B} \int_{[0,BN)} e^{-2\pi i\frac{\omega}{B} (\lambda - \sigma)} \sum_{j\in\mathbb Z} \widehat{\mathcal W}_{\left[ \frac{n}{2N} \right]}\left(\frac{\omega}{B}+Nj\right) \overline{\widehat{\mathcal W}_{\left[ \frac{n}{2N} \right]}\left(\frac{\omega}{B}+Nj\right)}\,d\omega\\
 & = \frac{e^{-i \pi \frac{A}{B}(\lambda^2 - \sigma^2)}}{B} \int_{[0,BN)} e^{-2\pi i\frac{\omega}{B} (\lambda - \sigma)} \sum_{j\in\mathbb Z} \left| \Lambda^M_k\left(\frac{\omega}{2NB} + \frac{j}{2} \right) \right|^2 \, \left| \widehat{\mathcal W}_{\left[ \frac{n}{2N} \right]} \left(\frac{\omega}{2NB} + \frac{j}{2}\right) \right|^2d\omega \\
 & = \frac{e^{-i \pi \frac{A}{B}(\lambda^2 - \sigma^2)}}{B} \int_{[0,BN)} e^{-2\pi i\frac{\omega}{B} (\lambda - \sigma)}  \left|\Lambda^M_k \left(\frac{\omega}{2NB}\right) \right|^2 \sum_{j\in\mathbb Z} \left| \widehat{\mathcal W}_{\left[ \frac{n}{2N} \right]}\left(\frac{\omega}{2NB} + \frac{j}{2}\right) \right|^2 d\omega\\
 &=\frac{e^{-i \pi \frac{A}{B}(\lambda^2 - \sigma^2)}}{B}\int_{[0,BN)} e^{-i2\pi \frac{\omega}{B} (\lambda - \sigma)}  \left|\Lambda^M_k \left( \frac{\omega}{2NB}\right) \right|^2 h_n \left(\frac{\omega}{2NB}\right)d\omega,
\end{align*}
where $h_n(\omega)=\sum_{j\in\mathbb Z} \big| \widehat{\mathcal W}_{\left[ \frac{n}{2N} \right]}\left(\omega+ Nj \right) \big|^2.$ In view of specific form of $\Omega$, we can write
\begin{align*}
\left| \Lambda^M_k \Big( \frac{\omega}{2NB} \Big) \right|^2 & =\left\{\Lambda^{M,1}_k\Big( \frac{\omega}{2NB} \Big) + e^{-i2\pi \frac{\omega}{B} \frac{r}{N}}\Lambda^{M,2}_k\Big( \frac{\omega}{2NB} \Big) \right\} \left\{ \overline{\Lambda^{M,1}_k\Big( \frac{\omega}{2NB} \Big)} + e^{i2\pi \frac{\omega}{B} \frac{r}{N}}  \overline{\Lambda^{M,2}_k\Big( \frac{\omega}{2NB} \Big)} \right\}  \\
  &=\left|\Lambda^{M,1}_k\Big( \frac{\omega}{2NB} \Big) \right|^2 + \left| \Lambda^{M,2}_k\Big( \frac{\omega}{2NB} \Big) \right|^2 + \overline{\Lambda^{M,1}_k\Big( \frac{\omega}{2NB} \Big) } \, \Lambda^{M,2}_k\Big( \frac{\omega}{2NB} \Big) \, e^{-i2\pi \frac{\omega}{B} \frac{r}{N}}\\
   &\qquad+ \Lambda^{M,1}_k \Big( \frac{\omega}{2NB} \Big) \, \overline{\Lambda^{M,2}_k\Big( \frac{\omega}{2NB} \Big)} \, e^{i2\pi \frac{\omega}{B} \frac{r}{N}}.
\end{align*}
Using (4.6) for $ \lambda=2m_{1},\sigma=2m_{2}, m_1,m_2\in\mathbb Z$, we obtain
\begin{align*}
&\Big\langle \mathcal W_n (t-\lambda),\mathcal W_n (t-\sigma) \Big\rangle \\
&=\frac{e^{-i 4\pi \frac{A}{B}( m^2_1 - m^2_2)}}{B} \Bigg[ \int_{[0,BN)} e^{-i4\pi \frac{\omega}{B} (m_1 - m_2)} \left\{ \left|\Lambda^{M,1}_k\left( \frac{\omega}{2NB}\right)\right|^2 + \left|\Lambda^{M,2}_k\left(\frac{\omega}{2NB}\right) \right|^2   \right\} h_n \left(\frac{\omega}{2NB} \right)d\omega\\
&\quad+ \int_{[0,BN)} e^{-i4\pi \frac{\omega}{B} (m_1 - m_2)} \, \overline{\Lambda^{M,1}_k\left(\frac{\omega}{2NB} \right) } \, \Lambda^{M,2}_k\left(\frac{\omega}{2NB}\right) \, e^{-i2\pi \frac{\omega}{B} \frac{r}{N}} \, h_n \left(\frac{\omega}{2NB}\right) d\omega\\
& \quad + \int_{[0,BN)} e^{-i4\pi \frac{\omega}{B} (m_1 - m_2)} \, \Lambda^{M,1}_k \left(\frac{\omega}{2NB} \right) \, \overline{\Lambda^{M,2}_k\left(\frac{\omega}{2NB} \right)} \, e^{i2\pi \frac{\omega}{B} \frac{r}{N}} \,  h_n \left(\frac{\omega}{2NB} \right)d\omega\Bigg]\\
&= \frac{e^{-i \pi \frac{A}{B}(\lambda^2 - \sigma^2)}}{B}
\end{align*}
\begin{align*}
 &\quad\times \Bigg[\int_{[0,\frac{B}{2})} e^{-i4\pi \frac{\omega}{B} (m_1 - m_2)} \left\{ \left|\Lambda^{M,1}_k\left( \frac{\omega}{2NB}\right)\right|^2 + \left|\Lambda^{M,2}_k\left(\frac{\omega}{2NB}\right) \right|^2   \right\} \sum_{p=0}^{2N-1} h_n \left(\frac{\omega}{2NB}+ \frac{p}{4N} \right)d\omega\\
 &\quad+ \int_{[0,\frac{B}{2})} e^{-i4\pi \frac{\omega}{B} (m_1 - m_2)}\overline{\Lambda^{M,1}_k\left(\frac{\omega}{2NB}\right) } \, \Lambda^{M,2}_k\left(\frac{\omega}{2NB}\right) \, e^{-i2\pi \frac{\omega}{B} \frac{r}{N}} \sum_{p=0}^{2N-1} e^{-i\pi\frac{r}{N} p} \, h_n \Big( \frac{\omega}{2NB} + \frac{p}{4N} \Big) \,d\omega\\
 &\quad+ \int_{[0,\frac{B}{2})} e^{-i4\pi \frac{\omega}{B} (m_1 - m_2)} \, m^{M,1}_\ell \Big( \frac{\omega}{2NB} \Big) \, \overline{ m^{M,2}_\ell \Big( \frac{\omega}{2NB} \Big)} \, e^{i2\pi \frac{\omega}{B} \frac{r}{N}}\sum_{p=0}^{2N-1} e^{-i\pi\frac{r}{N} p} \, h_n \Big( \frac{\omega}{2NB} + \frac{p}{4N} \Big) \, d\omega\Bigg] \\
&  = 2\,e^{-i 4\pi \frac{A}{B}( m^2_1 - m^2_2)} \, \int_{[0,\frac{B}{2})} e^{-i4\pi \frac{\omega}{B} (m_1 - m_2)} \,d\omega\\
&  =  \delta_{m_1,m_2} \\
& = \delta_{\lambda,\sigma}.
\end{align*}
Similarly, for $ \lambda=2m_{1}+r/N,\sigma=2m_{2}$, where $m_1,m_2\in\mathbb Z$, again using (4.6), we have
\begin{align*}
 &\Big\langle \mathcal W_n (t- \lambda), \mathcal W_n (t- \sigma) \Big\rangle\\
&= \frac{e^{-i\pi \frac{A}{B} \left( (2m_1+r/N)^2 - (2m_2)^2 \right)}}{B}\\
 & \times\Bigg[ \int_{[0,\frac{B}{2})} e^{-i4\pi \frac{\omega}{B} (m_1 - m_2)} \left\{ \left|\Lambda^{M,1}_k \left(\frac{\omega}{2NB}\right)\right|^2 + \left| \Lambda^{M,2}_k\left(\frac{\omega}{2NB}\right) \right|^2 \right\}\\
 &\qquad\times e^{i2\pi r\omega/NB }\sum_{p=0}^{2N-1} e^{-i\pi rp/N} \, h_n \left(\frac{\omega}{2NB} + \frac{p}{4N}\right)d\omega\\
 &+ \int_{[0,\frac{B}{2})} e^{-i4\pi \frac{\omega}{B} (m_1 - m_2)}\,\overline{\Lambda^{M,1}_k \left(\frac{\omega}{2NB}\right) }\, \Lambda^{M,2}_k\left(\frac{\omega}{2NB} \right)\sum_{p=0}^{2N-1} e^{-i\pi rp/N} \, h_n \left(\frac{\omega}{2NB} + \frac{p}{4N}\right)d\omega\\
  &+ \int_{[0,\frac{B}{2})} e^{-i4\pi \frac{\omega}{B} (m_1 - m_2)}\Lambda^{M,1}_k\left(\frac{\omega}{2NB}\right) \overline{ \Lambda^{M,2}_k\left(\frac{\omega}{2NB} \right)} \, e^{i4\pi r\omega/NB } \sum_{p=0}^{2N-1} e^{-i\pi rp/N}\, h_n \left(\frac{\omega}{2NB} + \frac{p}{4N}\right)d\omega\Bigg] \\
& = 0.
\end{align*}

This completes the proof of the Theorem 2.2. \quad\fbox

\parindent=0mm\vspace{.0in}

{\bf Theorem 4.3.} {\it For any $n_1,n_2\in\mathbb Z^+$ and $\lambda ,\sigma \in\Lambda$, we have}
\begin{align*}
  \Big\langle \mathcal W_{n_1} (t- \lambda),\mathcal W_{n_2} (t- \sigma) \Big\rangle = \delta_{n_1,n_2}\delta_{\lambda,\sigma},\tag{4.10}
\end{align*}
{\it  where $\big\{\mathcal W_n : n \ge 0\big\}$ is the  nonuniform wavelet packet associated with linear canonical transform.}

\parindent=0mm \vspace{.1in}
{\it Proof.} For $n_1 = n_2$, the result follows immediately from  Theorem 4.2. For $n_1 \neq n_2$ with  $n_1\geq n_2$, we assume that $n_1 = 2N \left[\frac{n_1}{2N} \right]+k,$ and  $n_2 = 2N \left[ \frac{n_2}{2N} \right]+\ell$, where $0\leq k,\ell \leq 2N-1$.

\parindent=0mm \vspace{.1in}
{\it Case-I.} If $\left[ \frac{n_1}{2N} \right] = \left[ \frac{n_2}{2N} \right]$, then for $k\neq \ell$, we have
\begin{align*}
&\Big\langle \mathcal W_{n_1} (t- \lambda),\mathcal W_{n_2} (t- \sigma) \Big\rangle =\\
 &= \frac{e^{-i \pi \frac{A}{B}( \lambda^2 - \sigma^2)}}{B} \int_{\mathbb R} e^{-i2\pi \frac{\omega}{B} (\lambda - \sigma)} \,\widehat{\mathcal W}_{n_1}\left(\frac{\omega}{B}\right)\overline{\widehat{\mathcal W}_{n_2}\left(\frac{\omega}{B}\right)} \,d\omega \\
 &= \frac{e^{-i \pi \frac{A}{B}( \lambda^2 - \sigma^2)}}{B} \int_{[0,BN)} e^{-i2\pi i \frac{\omega}{B} (\lambda - \sigma)}  \sum_{j\in\mathbb Z} \widehat{\mathcal W}_{n_1}\left(\frac{\omega}{B}+Nj\right)\overline{\widehat{\mathcal W}_{n_1}\left(\frac{\omega}{B}+Nj\right)} \,d\omega
 \end{align*}
 \begin{align*}
 &= \frac{e^{-i \pi \frac{A}{B}( \lambda^2 - \sigma^2)}}{B}  \int_{[0,BN)} e^{-i2\pi \frac{\omega}{B} (\lambda - \sigma)} h\left(\frac{\omega}{B}\right)d\omega,
\end{align*}
where $h(\omega) = \sum_{j\in\mathbb Z}\widehat{\mathcal W}_{n_1}(\omega+Nj)\overline{\widehat{\mathcal W}_{n_1}(\omega+Nj)}$. Consequently, we have
\begin{align*}
   h\left(\frac{2N\omega}{B}\right)&= \sum_{j\in\mathbb Z}\widehat{\mathcal W}_{n_1}\left(2N\Big(\frac{\omega}{B}+\frac{j}{2}\Big)\right)\overline{\widehat{\mathcal W}_{n_1}\left(2N\Big(\frac{\omega}{B}+\frac{j}{2}\Big)\right)}\\
   & = \sum_{j\in\mathbb Z} \Lambda_k^M \left(\frac{\omega}{B} + \frac{j}{2} \right) \widehat{\mathcal W}_{\left[ \frac{n_1}{2N} \right]} \left(\frac{\omega}{B} + \frac{j}{2} \right) \, \overline{\Lambda_\ell^M \left(\frac{\omega}{B} + \frac{j}{2}\right) \widehat{\mathcal W}_{\left[ \frac{n_2}{2N} \right]} \left( \frac{\omega}{B} + \frac{j}{2} \right)}   \\
   & = 2 \sum_{j\in\mathbb Z}  \Lambda_k^M \left(\frac{\omega}{B} + \frac{j}{2}\right) \, \overline{\Lambda_\ell^M \left(\frac{\omega}{B} + \frac{j}{2} \right) }.
\end{align*}
Therefore,
\begin{align*}
\sum_{p=0}^{2N-1}h\left(\frac{\omega}{B}+ \frac{p}{2}\right) &= \sum_{j\in\mathbb Z} \widehat{\mathcal W}_{n_1}\left(\frac{\omega}{B}+ \frac{j}{2}\right) \, \overline{\widehat{\mathcal W}_{n_2}\left(\frac{\omega}{B}+ \frac{j}{2}\right)} \\
& = 2 \sum_{j=0}^{2N-1} \Lambda_k^M \left(\frac{\omega}{2NB} + \frac{j}{4N} \right) \, \overline{\Lambda_\ell^M \left(\frac{\omega}{2NB} + \frac{j}{4N}\right) }.
\end{align*}
Assume that $\lambda =2m_1 $ and $ \sigma =2m_2$, where $m_1,m_2\in \mathbb Z$, then we have
\begin{align*}
& \Big\langle\mathcal W_{n_1} (t- \lambda), \mathcal W_{n_2} (t- \sigma) \Big\rangle\\
&= \frac{e^{-i 4\pi \frac{A}{B}( m_1^2 - m_2^2 )}}{B} \int_{[0,BN)} e^{-i4\pi \frac{\omega}{B} (m_1 - m_2)} \, h\left(\frac{\omega}{B}\right)d\omega\\
&= \frac{e^{-i 4\pi \frac{A}{B}( m_1^2 - m_2^2 )}}{B}\int_{[0,\frac{B}{2})} e^{-i4\pi \frac{\omega}{B} (m_1 - m_2)} \left[ \sum_{p=0}^{2N-1} h \left(\frac{\omega}{B} + \frac{p}{2} \right) \right]d\omega\\
&=\frac{2\,e^{-i 4\pi \frac{A}{B}( m_1^2 - m_2^2 )}}{B}\int_{[0,\frac{B}{2})} e^{-i4\pi \frac{\omega}{B} (m_1 - m_2)} \left[ \sum_{j=0}^{2N-1} \Lambda_k^M \left( \frac{\omega}{2NB} + \frac{j}{4N}\right)\overline{\Lambda_\ell^M \left( \frac{\omega}{2NB} + \frac{j}{4N} \right) } \right]d\omega\\
 &= \delta_{m_1,m_2} \, \delta_{k,\ell}  \\
 &= \delta_{\lambda,\sigma} \, \delta_{k,\ell}.
\end{align*}
Similarly, if $\lambda = 2m_1 + r/N$ and $ \sigma = 2m_2$, where $m_1,m_2 \in \mathbb Z$, then by virtue of (4.6), we have
\begin{align*}
& \Big\langle\mathcal W_{n_1} (t- \lambda), \mathcal W_{n_2} (t- \sigma) \Big\rangle\\
&= \frac{e^{-i\pi\frac{A}{B} \left( (2m_1+r/N)^2 - (2m_2)^2 \right)}}{B} \int_{[0,BN)} e^{-i4\pi \frac{\omega}{B} (m_1 - m_2)} \, e^{-i2\pi \frac{\omega}{B} \frac{r}{N}}\, h\left(\frac{\omega}{B}\right)d\omega  \\
&= \frac{e^{-i\pi\frac{A}{B} \left( (2m_1+r/N)^2 - (2m_2)^2 \right)}}{B}\int_{[0,\frac{B}{2})} e^{-i4\pi \frac{\omega}{B} (m_1 - m_2)}  \, e^{-i2\pi \frac{\omega}{B} \frac{r}{N}} \left[ \sum_{p=0}^{2N-1} e^{-i\pi rp/N} h\left(\frac{\omega}{B} + \frac{p}{2} \right) \right] d\omega\\
&= 0.
\end{align*}

{\it Case-II.} In case $\left[ \frac{n_1}{2N}\right] = \left[ \frac{n_2}{2N}\right]$, then  we can take $\left[ \frac{n_1}{2N}\right] = 2N \left[ \frac{\left[  n_1/2N  \right]}{2N} \right] + k_1, \left[ \frac{n_2}{2N}\right] = 2N \left[ \frac{\left[  n_2/2N  \right]}{2N} \right] + \ell_1, ~0\leq k_1,\ell_1 \leq 2N-1.$ Let $\left[ \frac{n_1}{2N}\right] = 2Np_1 + k_1$ and $\left[ \frac{n_2}{2N}\right] = 2Nq_1 + \ell_1$, where $p_1 = \left[ \frac{\left[  n_1/2N  \right]}{2N} \right] $ and $q_1 = \left[ \frac{\left[  n_2/2N  \right]}{2N} \right]$. In case $p_1 = q_1$, then the result follows from  Case-I. If $p_1 \neq q_1$, then we take $p_1 = 2N \left[ \frac{p_1}{2N} + k_2 \right] = 2Np_2 +k_2 $ and $q_1 = 2N \left[ \frac{q_1}{2N} + \ell_2 \right] = 2Nq_2 +\ell_2$, where $ 0\leq k_2,\ell_2 \leq 2N-1$. Again, if $p_2 = q_2$, then the result follows from Case-I. However, if $p_2\neq q_2$, then apply the above procedure. After performing a finite number of steps, we obtain
$p_{m-1} = 2Np_m + k_m$ and $q_{m-1} = 2Nq_m + \ell_m$, where $0\leq p_m,q_m \leq 2N-1$. Again, there arises two case:  when $p_m = q_m$, then there is nothing to prove as the result follows again by applying Case-I. But, in case $p_m \neq q_m$, then we have
\begin{align*}
 &\Big\langle \mathcal W_{n_1}(t-\lambda) , \mathcal W_{n_2} (t-\sigma) \Big\rangle \\
 &= \frac{e^{-i \pi \frac{A}{B}( \lambda^2 - \sigma^2)}}{B} \int_{\mathbb R} e^{-i2\pi \frac{\omega}{B} (\lambda - \sigma)} \, \widehat{\mathcal W}_{n_1}\left(\frac{\omega}{B}\right)\overline{ \widehat{\mathcal W}_{n_2}\left(\frac{\omega}{B}\right)}\,d\omega\\
 &=\frac{e^{-i \pi \frac{A}{B}( \lambda^2 - \sigma^2)}}{B} \int_{\mathbb R} e^{-i2\pi \frac{\omega}{B} (\lambda - \sigma)} \Lambda^M_{k_1} \Big( \frac{\omega}{2NB}\Big)\widehat{\mathcal W}_{\left[ \frac{n_1}{2N}\right]} \Big( \frac{\omega}{2NB}\Big)\overline{\Lambda^M_{\ell_1} \Big( \frac{\omega}{2NB}\Big) \, \widehat{\mathcal W}_{\left[ \frac{n_2}{2N}\right]} \Big( \frac{\omega}{2NB}\Big)}\,d\omega\\
 &= \frac{e^{-i \pi \frac{A}{B}( \lambda^2 - \sigma^2)}}{B} \int_{\mathbb R} e^{-i2\pi \frac{\omega}{B} (\lambda - \sigma)} \Lambda^M_{k_1} \left( \frac{\omega}{2NB}\right) \Lambda^M_{k_2} \left( \frac{\omega}{(2NB)^2}\right) \dots  \Lambda^M_{k_m}\left( \frac{\omega}{(2NB)^m}\right) \\
 &\quad \widehat{\mathcal W}_{p_m} \left( \frac{\omega}{(2NB)^m} \right) \,\overline{\Lambda^M_{\ell_1} \left( \frac{\omega}{2NB}\right)} \, \overline{\Lambda^M_{\ell_2} \left( \frac{\omega}{(2NB)^2}\right)} \dots \overline{\Lambda^M_{\ell_m} \left( \frac{\omega}{(2NB)^m}\right) } \overline{ \widehat{\mathcal W}_{q_m} \left( \frac{\omega}{(2NB)^m} \right)} d\omega\\
 &=\frac{e^{-i \pi \frac{A}{B}( \lambda^2 - \sigma^2)}}{B} \int_{\mathbb R} e^{-i2\pi \frac{\omega}{B} (\lambda - \sigma)}d\omega \\
  &\quad\left[ \prod_{n=1}^{m} \Lambda^M_{k_n} \left( \frac{\omega}{(2NB)^n}\right) \right] \widehat{\mathcal W}_{p_m} \left( \frac{\omega}{(2NB)^m} \right)\left[ \prod_{n=1}^{m} \overline{ \Lambda^M_{\ell_n} \left( \frac{\omega}{(2NB)^n}\right)} \right] \overline{ \widehat{\mathcal W}_{q_m} \left( \frac{\omega}{(2NB)^m} \right)} d\omega  \\
 &= \frac{e^{-i \pi \frac{A}{B}( \lambda^2 - \sigma^2)}}{B} \int_{[0,BN)} e^{-i2\pi \frac{\omega}{B} (\lambda - \sigma)}\,d\omega\\
 &\left[ \prod_{n=1}^{m} \Lambda^M_{k_n} \left( \frac{\omega}{(2NB)^n}\right) \right] \left[ \sum_{j\in\mathbb Z} \widehat{\mathcal W}_{p_m} \left( \frac{\omega}{(2NB)^m} + Nj \right) \overline{ \widehat{\mathcal W}_{q_m} \left( \frac{\omega}{(2NB)^m} + Nj \right)} \right] \left[ \prod_{n=1}^{m} \overline{ \Lambda^M_{\ell_n} \left( \frac{\omega}{(2NB)^n}\right) }\right] \\
 &= 0.
\end{align*}

This completes the proof the Theorem 4.3.\quad\fbox

\parindent=0mm\vspace{.1in}

{\bf Theorem 4.4.} {\it   Let $\big\{\mathcal W_n:n\ge 0 \big\}$ be the basic nonuniform wavelet packet associated with an LCT-NUMRA $\big\{V_{j}^{M}: j\in\mathbb Z\big\}$. Then, for all $n\in \mathbb Z^+$ and $k,\ell\in\big\{0,1,2,\dots, 2N-1\big\}$, we have}
\begin{align*}
  \Big\langle \mathcal W_{2Nn+k}(t-\lambda),\mathcal W_{2Nn+\ell}(t-\sigma)\Big\rangle=\delta_{k,\ell}\,\delta_{\lambda,\sigma},~~\lambda,\sigma\in\Omega.\tag{4.11}
\end{align*}

\parindent=0mm \vspace{.1in}
{\it Proof.} Since $\{\psi^M_{k,0,\lambda}:1\leq k\leq2N-1,\lambda\in\Omega\}$ constitutes an orthonormal basis for $W_0^M$, so we have
\begin{align*}
&\Big\langle \mathcal W_{2Nn+k} (t- \lambda),\mathcal W_{2Nn+\ell} (t-\sigma) \Big\rangle \\
 &= \frac{e^{-i\pi\frac{A}{B}\left(\lambda^2-\sigma^2\right)}}{B}\int_{\mathbb R}e^{-i2\pi \frac{\omega}{B}(\lambda-\sigma)}\,\widehat{\mathcal W}_{2Nn+k}\left(\frac{\omega}{B}\right)\overline{\widehat{\mathcal W}_{2Nn+\ell}\left(\frac{\omega}{B}\right)}\,d\omega  \\
 &= \frac{e^{-i \pi \frac{A}{B}( \lambda^2 - \sigma^2)}}{B} \int_{[0,BN)}e^{-i2\pi\frac{\omega}{B}(\lambda-\sigma)}\sum_{j\in\mathbb Z}\widehat{\mathcal W}_{2Nn+k}\left(\frac{\omega}{B}+Nj\right)\overline{\widehat{\mathcal W}_{2Nn+\ell}\left(\frac{\omega}{B}+Nj\right)}\,d\omega\\
 &=\frac{e^{-i \pi \frac{A}{B}( \lambda^2 - \sigma^2)}}{B} \int_{[0,BN)}e^{-i2\pi\frac{\omega}{B}(\lambda-\sigma)}h\left(\frac{\omega}{B}\right)d\omega,
\end{align*}
where $h(\omega)= \sum_{j\in\mathbb Z}\widehat{\mathcal W}_{2Nn+k}(\omega+Nj)\,\overline{\widehat{\mathcal W}_{2Nn+\ell}(\omega+Nj)}$. By continuing the same way as in Theorem 4.3, we may obtain

\begin{align*}
   h\left(\frac{2N\omega}{B}\right)=2\sum_{j\in\mathbb Z}\Lambda_k^M\left(\frac{\omega}{B}+\frac{j}{2}\right)\,\overline{\Lambda_\ell^M\left(\frac{\omega}{B}+\frac{j}{2}\right) }.
\end{align*}
Therefore, we have
\begin{align*}
\sum_{p=0}^{2N-1}h\left(\frac{\omega}{B}+\frac{p}{2}\right)=2\sum_{j=0}^{2N-1}\Lambda_k^M\left(\frac{\omega}{2NB}+\frac{j}{4N}\right)\, \overline{\Lambda_\ell^M\left(\frac{\omega}{2NB}+\frac{j}{4N}\right)}.
\end{align*}
Assume that $\lambda=2m_1$ and $\sigma=2m_2$, where $m_1,m_2\in \mathbb Z$, then we have
\begin{align*}
& \Big\langle\mathcal W_{n_1} (t- \lambda), \mathcal W_{n_2} (t- \sigma) \Big\rangle\\
&= \frac{e^{-i4\pi\frac{A}{B}(m_1^2-m_2^2)}}{B}\int_{[0,BN)}e^{-i4\pi\frac{\omega}{B}(m_1-m_2)}\,h\left(\frac{\omega}{B}\right)d\omega\\
&= \frac{e^{-i4\pi\frac{A}{B}(m_1^2-m_2^2)}}{B}\int_{[0,\frac{B}{2})}e^{-i4\pi\frac{\omega}{B}(m_1-m_2)}\left[\sum_{p=0}^{2N-1} h \left(\frac{\omega}{B}+\frac{p}{2}\right)\right]d\omega\\
&= \frac{2\,e^{-i4\pi\frac{A}{B}(m_1^2-m_2^2)}}{B}\int_{[0,\frac{B}{2})}e^{-i4\pi\frac{\omega}{B}(m_1-m_2)}\left[\sum_{j=0}^{2N-1} \Lambda_k^M\left(\frac{\omega}{2NB}+\frac{j}{4N}\right)\overline{\Lambda_\ell^M\left(\frac{\omega}{2NB}+\frac{j}{4N}\right)} \right]d\omega\\
 &= \delta_{m_1,m_2}\,\delta_{k,\ell}  \\
 &= \delta_{\lambda,\sigma}\,\delta_{k,\ell}.
\end{align*}
Similarly, if $\lambda=2m_1+r/N$ and $\sigma=2m_2$, where $m_1,m_2\in\mathbb Z$, then by virtue of (4.6), we have
\begin{align*}
& \Big\langle\mathcal W_{n_1}(t-\lambda),\mathcal W_{n_2}(t-\sigma)\Big\rangle\\
&= \frac{e^{-i\pi\frac{A}{B}\left((2m_1+r/N)^2-(2m_2)^2\right)}}{B}\int_{[0,BN)}e^{-i4\pi\frac{\omega}{B}(m_1-m_2)}\,e^{-i2\pi\frac{\omega}{B} \frac{r}{N}}\,h\left(\frac{\omega}{B}\right)d\omega\\
&= \frac{e^{-i\pi\frac{A}{B}\left((2m_1+r/N)^2-(2m_2)^2\right)}}{B}\int_{[0,\frac{B}{2})}e^{-i4\pi\frac{\omega}{B}(m_1-m_2)}\,e^{-i2\pi\frac{\omega}{B} \frac{r}{N}}\left[\sum_{p=0}^{2N-1}e^{-i\pi p\frac{r}{N}}h\Big(\frac{\omega}{B}+\frac{p}{2}\Big)\right]d\omega\\
&= 0.
\end{align*}
Thus, the system  $\big\{\mathcal W_n (t - \lambda): n \in \mathbb Z^+ \big\} $ constitutes an orthonormal basis in $L^2(\mathbb R)$. \quad\fbox

\section{Orthogonal Decomposition of LCT Nonuniform Wavelet Packet Subspaces}

\parindent=0mm\vspace{.0in}
Let $\{\mathcal W_n:n\ge0\}$ be the family of nonuniform wavelet packets associated with  LCT-NUMRA $\big\{V_j^M:j\in\mathbb Z\big\}$ of $L^2(\mathbb R)$. For $n\ge0$ and $j\in\mathbb Z$, we consider the family of subspaces
\begin{align*}
U_{j,n}^M=\overline{\text{span}}\left\{(2N)^{j/2}\mathcal W_n\big((2N)^jt-\lambda\big)\,e^{-i\pi\frac{A}{B}(t^2-\lambda^2)}:\lambda\in\Omega\right\}.\tag{5.1}
\end{align*}
When $n=0$ and $1\le n\le2N-1$ in (5.1), it is easy to see that
\begin{align*}
U_{j,0}^M=V_j^M\quad\text{and}\quad\bigoplus_{m=1}^{2N-1}U_{j,m}^M=\mathcal W_j^M,\quad j\in\mathbb Z.\tag{5.2}
\end{align*}
Thus, the orthogonal decomposition $V_j^M=V_{j-1}^M\oplus W_{j-1}^M$ can be expressed as
\begin{align*}
U_{j,0}^M=U_{j-1,0}^M\bigoplus_{m=1}^{2N-1}U_{j-1,m}^M=\bigoplus_{m=0}^{2N-1}U_{j-1,m}^M,\quad j\in\mathbb Z\tag{5.3}
\end{align*}
In fact, this orthogonal decomposition can be generalized to any LCT-wavelet packet space $U_{j,n}^M$ that is spanned by
an orthogonal basis of functions $\big\{ \mathcal W_{n,j,\lambda}^{M}(t): n\ge 0, j\in\mathbb Z, \lambda\in \Omega\big\}$.

\parindent=0mm\vspace{.1in}
{\bf Theorem 5.1.} {\it For  $n\ge0$ and $j\in\mathbb Z$, following relationship holds}
\begin{align*}
U_{j,n}^M=\bigoplus_{\ell=0}^{2N-1}U_{j,2nN+\ell}^M~.\tag{5.4}
\end{align*}

\parindent=0mm\vspace{.1in}
{\it Proof.}  From the   recursive property of the Definition 4.1, we observe that  $U_{j-1,2nN+\ell}^M,\;0\le\ell\le2N-1$ are subspaces of $U_{j,n}^M$. Thus, we have
\begin{align*}
\bigoplus_{\ell=0}^{2N-1}U_{j-1,2nN+\ell}^M\subseteq U_{j,n}^M,\quad j\in\mathbb Z.\tag{5.5}
\end{align*}
Theorem 4.3 implies that $U_{j-1,n_1}^M$ and $U_{j-1,n_2}^M$ are orthogonal to each other. So, it suffices to show that for any $f\in U_{j,n}^M$, if
\begin{align*}
\Big\langle f,\mathcal W_{j-1,\lambda,p}^M\Big\rangle=0,\quad p=2nN+k,\;0\le k\le2N-1,\tag{5.6}
\end{align*}
then $f(t)=0$. We now proceed to prove the result as follows:

\parindent=8mm\vspace{.1in}
For any $f(t)\in U_{j,n}^M$, we see that
\begin{align*}
f(t)=\sum_{\lambda\in\Omega}c_\lambda(2N)^{j/2}\,\mathcal W_n\big((2N)^jt-\lambda\big)\,e^{-i\pi\frac{A}{B}(t^2-\lambda^2)}.\tag{5.7}
\end{align*}
Applying LCT on both sides of equation (5.7) results in
\begin{align*}
\mathscr L_M\big[f\big](\omega)&=\int_{\mathbb R}\sum_{\lambda\in\Omega}c_\lambda(2N)^{j/2}\,\mathcal W_n\big((2N)^jt-\lambda\big)\,e^{-i\pi\frac{A}{B}(t^2-\lambda^2)}\,\mathcal K_M(\omega,t)\,dt\\
&=\frac{(2N)^{-j/2}}{\sqrt{iB}}\,e^{i\pi D\frac{\omega^2}{B}}\sum_{\lambda\in\Omega}c_\lambda\,e^{i\pi\frac{A}{B}\lambda^2}e^{-i2\pi\lambda\frac{\omega}{(2N)^jB}}\,\widehat{\mathcal W}_n\left(\frac{\omega}{(2N)^jB}\right)\\
&=\frac{(2N)^{-j/2}}{\sqrt{iB}}\,e^{i\pi D\frac{\omega^2}{B}}\sum_{\lambda\in\Omega}c_\lambda^M e^{-i2\pi\lambda\frac{\omega}{(2N)^jB}}\,\widehat{\mathcal W}_n\left(\frac{\omega}{(2N)^jB}\right)\\
&=\frac{(2N)^{-j/2}}{\sqrt{iB}}\,e^{i\pi D\frac{\omega^2}{B}}\,\widetilde{C}\left(\frac{\omega}{(2N)^jB}\right)\widehat{\mathcal W}_n\left(\frac{\omega}{(2N)^jB}\right),\tag{5.8}
\end{align*}
where $\widetilde{C}(\omega)=\sum_{\lambda\in\Omega}c_\lambda^M e^{-i2\pi\lambda\omega}$. Moreover, we have
\begin{align*}
\mathscr L_M\Big[\mathcal W_{j-1,\lambda,p}^M\Big](\omega)&=\int_{\mathbb R}(2N)^{(j-1)/2}\,\mathcal W_p\big((2N)^{j-1}t-\lambda\big)\,e^{-\pi\frac{A}{B}(t^2-\lambda^2)}\,\mathcal K_M(\omega,t)\,dt\\
&=\frac{(2N)^{-(j-1)/2}}{\sqrt{iB}}\,e^{i\frac{\pi}{B}\big(A\lambda^2-2\lambda\frac{\omega}{(2N)^{j-1}}+D\omega^2\big)}\,\widehat{\mathcal W}_p\left(\frac{\omega}{(2N)^{j-1}B}\right),\tag{5.9}
\end{align*}
where $p=2nN+k,\;0\le k\le2N-1$. Then, combining the relations (5.8), (5.9) and Parseval’s identity of the LCT, the L.H.S of (5.6) can be expressed as
\begin{align*}
&\Big\langle f,\mathcal W_{j-1,\lambda,p}^M\Big\rangle\\
&\quad=\frac{(2N)^{-j+\frac{1}{2}}\,e^{-i\pi\frac{A}{B}\lambda^2}}{B}\int_\mathbb Re^{i2\pi\frac{\lambda}{(2N)^{j-1}}\frac{\omega}{B}}\,\widetilde{C}\left(\frac{\omega}{(2N)^jB}\right)\widehat{\mathcal W}_n\left(\frac{\omega}{(2N)^jB}\right)\widehat{\mathcal W}_p\left(\frac{\omega}{(2N)^{j-1}B}\right)d\omega\\
&\quad=\frac{(2N)^{-\frac{1}{2}}\,e^{-i\pi\frac{A}{B}\lambda^2}}{B}\int_\mathbb Re^{i2\pi\lambda\frac{\omega}{B}}\,\widetilde{C}\left(\frac{\omega}{2NB}\right)\widehat{\mathcal W}_n\left(\frac{\omega}{2NB}\right)\widehat{\mathcal W}_p\left(\frac{\omega}{B}\right)d\omega.\tag{5.10}
\end{align*}
Since $p=2nN+k,$ with $0\le k\le2N-1$, so we have $\left[\frac{p}{2N}\right]=n$. Then, combining (5.10) and (4.3) yields
\begin{align*}
&\Big\langle f,\mathcal W_{j-1,\lambda,p}^M\Big\rangle\\
&\qquad=\frac{(2N)^{-\frac{1}{2}}\,e^{-i\pi\frac{A}{B}\lambda^2}}{B}\int_\mathbb Re^{i2\pi\lambda\frac{\omega}{B}}\,\widetilde{C}\left(\frac{\omega}{2NB}\right)\widehat{\mathcal W}_n\left(\frac{\omega}{2NB}\right)\overline{\Lambda_k^M\left(\frac{\omega}{2NB}\right)\widehat{\mathcal W}_n\left(\frac{\omega}{2NB}\right)}\,d\omega\\
&\qquad=\frac{(2N)^{-\frac{1}{2}}\,e^{-i\pi\frac{A}{B}\lambda^2}}{B}\int_\mathbb Re^{i2\pi\lambda\frac{\omega}{B}}\,\widetilde{C}\left(\frac{\omega}{2NB}\right)\overline{\Lambda_k^M\left(\frac{\omega}{2NB}\right)}\left|\widehat{\mathcal W}_n\left(\frac{\omega}{2NB}\right)\right|^2d\omega\\
&\qquad=\frac{(2N)^{-\frac{1}{2}}\,e^{-i\pi\frac{A}{B}\lambda^2}}{B}\int_{[0,BN)}e^{i2\pi\lambda\frac{\omega}{B}}\,\widetilde{C}\left(\frac{\omega}{2NB}\right)\overline{\Lambda_k^M\left(\frac{\omega}{2NB}\right)}\sum_{j\in\mathbb Z}\left|\widehat{\mathcal W}_n\left(\frac{\omega}{2NB}+\frac{j}{2}\right)\right|^2d\omega\\
&\qquad=\frac{(2N)^{-\frac{1}{2}}\,e^{-i\pi\frac{A}{B}\lambda^2}}{B}\int_{[0,BN)}e^{i2\pi\lambda\frac{\omega}{B}}\,\widetilde{C}\left(\frac{\omega}{2NB}\right)\overline{\Lambda_k^M\left(\frac{\omega}{2NB}\right)}\,h_n\left(\frac{\omega}{2NB}\right)d\omega,\tag{5.11}
\end{align*}
where $h_n(\omega)=\sum_{j\in\mathbb Z} \big| \widehat{\mathcal W}_n\left(\omega+ Nj \right) \big|^2$. In view of specific form of $\Omega$ and using (4.6), relation (5.11) can be rewritten as
\begin{align*}
&\Big\langle f,\mathcal W_{j-1,\lambda,p}^M\Big\rangle\\
&=\frac{(2N)^{-\frac{1}{2}}\,e^{-i\pi\frac{A}{B}\lambda^2}}{B}\\
&\times \Bigg[\int_{[0,\frac{B}{2})}e^{i2\pi\lambda\frac{\omega}{B}}\,\widetilde{C}\left(\frac{\omega}{2NB}\right)\left\{ \left|\Lambda^{M,1}_k\left( \frac{\omega}{2NB}\right)\right|^2 + \left|\Lambda^{M,2}_k\left(\frac{\omega}{2NB}\right) \right|^2   \right\} \sum_{p=0}^{2N-1} h_n \left(\frac{\omega}{2NB}+ \frac{p}{4N} \right)d\omega\\
&+ \int_{[0,\frac{B}{2})}e^{i2\pi\lambda\frac{\omega}{B}}\,\widetilde{C}\left(\frac{\omega}{2NB}\right)\overline{\Lambda^{M,1}_k\left(\frac{\omega}{2NB}\right) } \, \Lambda^{M,2}_k\left(\frac{\omega}{2NB}\right) \, e^{-i2\pi \frac{\omega}{B} \frac{r}{N}} \sum_{p=0}^{2N-1} e^{-i\pi\frac{r}{N} p} \, h_n \Big( \frac{\omega}{2NB} + \frac{p}{4N} \Big) \, d\omega\\
&+ \int_{[0,\frac{B}{2})}e^{i2\pi\lambda\frac{\omega}{B}}\,\widetilde{C}\left(\frac{\omega}{2NB}\right)\, \Lambda^{M,1}_\ell \Big( \frac{\omega}{2NB} \Big) \, \overline{ \Lambda^{M,2}_\ell \Big( \frac{\omega}{2NB} \Big)} \, e^{i2\pi \frac{\omega}{B} \frac{r}{N}}\sum_{p=0}^{2N-1} e^{-i\pi\frac{r}{N} p} \, h_n \Big( \frac{\omega}{2NB} + \frac{p}{4N} \Big) \, d\omega \Bigg] \\
&= \frac{2\,(2N)^{-\frac{1}{2}}\,e^{-i\pi\frac{A}{B}\lambda^2}}{B}\int_{[0,\frac{B}{2})}e^{i2\pi\lambda\frac{\omega}{B}}\,\widetilde{C}\left(\frac{\omega}{2NB}\right)d\omega
\end{align*}
which implies that (5.6) is equivalent to $\widetilde{C}\left(\frac{\omega}{2NB}\right)=0$. Therefore, we have $\widetilde{C}\left(\frac{\omega}{B}\right)=0$, i.e., $c_\lambda=0$, which alongside (5.7)  yields $f(t)=0$. This completes the proof of Theorem 5.1.\quad\fbox

\parindent=0mm\vspace{.1in}

{\bf Theorem 5.2.} {\it For $j\ge0$, we have }
\begin{align*}
W_j^M=\displaystyle\bigoplus_{\ell=1}^{2N-1}U_{j,\ell}^M=\bigoplus_{\ell=2N}^{(2N)^2-1}U_{j-1,\ell}^M=\dots=\bigoplus_{\ell=(2N)^m }^{(2N)^{m+1}-1}U_{j-m,\ell}^M=\dots=\bigoplus_{\ell=(2N)^j}^{(2N)^{j+1}-1}U_{0,\ell}^M.\tag{5.12}
\end{align*}

\parindent=0mm\vspace{.0in}
{\it Proof.} The proof is obtained by repeated application of the previous theorem.\quad\fbox

\parindent=8mm\vspace{.1in}

It follows from Theorem 5.2 that for each $j\in\mathbb Z$, we have the following wavelet packet decomposition of $L^2(\mathbb R)$
\begin{align*}
L^2(\mathbb R)=\bigoplus_{j\in\mathbb Z}W_j^M=\bigoplus_{j\in\mathbb Z}\left(\bigoplus_{\ell=(2N)^m}^{(2N)^{m+1}-1}U_{j-m,\ell}^M\right)
\end{align*}
and subsequently, we  can construct many orthonormal bases of $L^2(\mathbb R)$.

\parindent=0mm\vspace{.1in}
{\it Data Availability Statement:}

\parindent=8mm\vspace{.1in}
The authors declare that there is no  supplementary data associated with the manuscript.

\end{document}